\documentclass[11pt]{amsart}
\allowdisplaybreaks

\usepackage{fancyhdr}
\usepackage[titletoc]{appendix}
\usepackage[dvipsnames]{xcolor}
\usepackage{etoolbox}
\usepackage{abstract}

\patchcmd{\chapter}{\normalfont}{\normalfont\color{BrickRed}}{}{}
\patchcmd{\section}{\normalfont}{\normalfont\color{RubineRed}}{}{}
\patchcmd{\subsection}{\normalfont}{\normalfont\color{blue}}{}{}

\usepackage{hyperref}
\hypersetup{colorlinks,linkcolor={BrickRed},citecolor={blue},urlcolor={red}}

\usepackage{color,xcolor}
\colorlet{BLUE}{blue}
\colorlet{RED}{red}
\colorlet{GRAY}{gray}
\colorlet{BROWN}{brown}
\definecolor{OliveGreen}{rgb}{0,0.6,0}

\newcommand{\Complex}{\mathbb{C}}
\newcommand{\Real}{\mathbb{R}}

\DeclareMathOperator{\diam}{diam}

\DeclareMathOperator{\Avg}{Avg}
\DeclareMathOperator{\Thin}{Thin}

\DeclareMathOperator{\dist}{dist}
\DeclareMathOperator{\Conv}{Conv}

\DeclareMathOperator{\Lip}{Lip}
\DeclareMathOperator{\CAT}{CAT}

\newcommand{\sr}{\stackrel}

\newcommand {\ol}{\overline}
\newcommand {\ul}{\underline}

\newcommand {\ra}{\rightarrow}
\newcommand {\Ra}{\Rightarrow}

\newcommand {\La}{\Leftarrow}

\newcommand {\ral}{\longrightarrow}

\newcommand {\vep}{\varepsilon}
\newcommand {\vphi}{\varphi}

\newcommand {\ld}{\lambda}

\newcommand {\al}{\alpha}

\newcommand{\wt}{\widetilde}

\newcommand {\bfr}{\begin{flushright}}
\newcommand {\efr}{\end{flushright}}
\newcommand {\bfl}{\begin{flushleft}}
\newcommand {\efl}{\end{flushleft}}

\newcommand {\nn} {\nonumber}

\newcommand {\txt}{\textrm}
\newcommand {\bd}{\begin{document}}
\newcommand {\ed}{\end{document}}

\newcommand {\bt}{\begin{tikzcd}}
\newcommand {\et}{\end{tikzcd}}
\newcommand {\be}{\begin{equation}}
\newcommand {\ee}{\end{equation}}
\newcommand {\bea}{\begin{eqnarray}}
\newcommand {\eea}{\end{eqnarray}}
\newcommand {\ba}{\begin{aligned}}
\newcommand {\ea}{\end{aligned}}

\newcommand{\bbib}{\begin{thebibliography}{99}}
\newcommand{\ebib}{\end{thebibliography}}
\newcommand {\bab}{\begin{abstract}}
\newcommand {\eab}{\end{abstract}}

\newcommand {\bc}{\begin{center}}
\newcommand {\ec}{\end{center}}
\newcommand {\bit}{\begin{itemize}}
\newcommand {\eit}{\end{itemize}}
\newcommand {\txtc}{\textcolor}

\def\A{\mathcal{A}}\def\B{\mathcal{B}}\def\C{\mathcal{C}}\def\D{\mathcal{D}}\def\E{\mathcal{E}}\def\F{\mathcal{F}}
\def\G{\mathcal{G}}\def\H{\mathcal{H}}\def\I{\mathcal{I}}\def\J{\mathcal{J}}\def\K{\mathcal{K}}\def\L{\mathcal{L}}
\def\M{\mathcal{M}}\def\N{\mathcal{N}}\def\O{\mathcal{O}}\def\P{\mathcal{P}}\def\Q{\mathcal{Q}}\def\R{\mathcal{R}}
\def\S{\mathcal{S}}\def\T{\mathcal{T}}\def\U{\mathcal{U}}\def\V{\mathcal{V}}\def\W{\mathcal{W}}\def\X{\mathcal{X}}
\def\Y{\mathcal{Y}}\def\Z{\mathcal{Z}}
\def\abar{\bar{a}}\def\bbar{\bar{b}}\def\cbar{\bar{c}}\def\dbar{\bar{d}}\def\ebar{\bar{e}}\def\fbar{\bar{f}}
\def\gbar{\bar{g}}\def\ibar{\bar{i}}\def\jbar{\bar{j}}\def\kbar{\bar{k}}\def\lbar{\bar{l}}
\def\mbar{\bar{m}}\def\nbar{\bar{n}}\def\obar{\bar{o}}\def\pbar{\bar{p}}\def\qbar{\bar{q}}\def\rbar{\bar{r}}
\def\sbar{\bar{s}}\def\tbar{\bar{t}}\def\ubar{\bar{u}}\def\vbar{\bar{v}}\def\wbar{\bar{w}}
\def\xbar{\bar{x}}\def\ybar{\bar{y}}\def\zbar{\bar{z}}
\def\Abar{\bar{A}}\def\Bbar{\bar{B}}\def\Cbar{\bar{C}}\def\Dbar{\bar{D}}\def\Ebar{\bar{E}}\def\Fbar{\bar{F}}
\def\Gbar{\bar{G}}\def\Hbar{\bar{H}}\def\Ibar{\bar{I}}\def\Jbar{\bar{J}}\def\Kbar{\bar{K}}\def\Lbar{\bar{L}}
\def\Mbar{\bar{M}}\def\Nbar{\bar{N}}\def\Obar{\bar{O}}\def\Pbar{\bar{P}}\def\Qbar{\bar{Q}}\def\Rbar{\bar{R}}
\def\Sbar{\bar{S}}\def\Tbar{\bar{T}}\def\Ubar{\bar{U}}\def\Vbar{\bar{V}}\def\Wbar{\bar{W}}
\def\Xbar{\bar{X}}\def\Ybar{\bar{Y}}\def\Zbar{\bar{Z}}
\def\0bar{\bar{0}}\def\1bar{\bar{1}}\def\2bar{\bar{2}}\def\3bar{\bar{3}}\def\4bar{\bar{4}}\def\5bar{\bar{5}}
\def\6bar{\bar{6}}\def\7bar{\bar{7}}\def\8bar{\bar{8}}\def\9bar{\bar{9}}
\def\cprod{\mathop{\prod\hspace{-0.457cm}\times}}
\def\tprod{\mathop{\prod\hspace{-0.457cm}\otimes}}
\def\dsum{\mathop{\sum\hspace{-0.458cm}\oplus}}
\def\Tr{\txt{Tr}}\def\tr{\txt{tr}}\def\str{\txt{Str}}
\def\bde{\begin{description}}\def\ede{\end{description}}

\renewcommand{\baselinestretch}{1.05}

\calclayout

\usepackage{chngcntr}
\counterwithout{equation}{section} 

\newtheorem{thm}{Theorem}[section]
\newtheorem{lmm}[thm]{Lemma}
\newtheorem{dfn}[thm]{Definition}
\newtheorem{crl}[thm]{Corollary}
\newtheorem{prp}[thm]{Proposition}
\newtheorem{example}[thm]{Example}
\newtheorem{rmk}[thm]{Remark}
\newtheorem{problem}[thm]{Problem}
\newtheorem{algorithm}[thm]{Algorithm}
\newtheorem{question}[thm]{Question}
\newtheorem{recall}[thm]{Recall}
\newtheorem{exercise}[thm]{Exercise}
\newtheorem{note}[thm]{Note}
\newtheorem{notation}[thm]{Notation}
\newtheorem{caution}[thm]{Caution}
\newtheorem{claim}[thm]{Claim}
\newtheorem{fact}[thm]{Fact}
\newtheorem{solution}[thm]{Solution}
\newtheorem{conjecture}[thm]{Conjecture}

\theoremstyle{definition}
\newtheorem{dfns}{Definitions}
\newtheorem{examples}{Examples}
\newtheorem{rmks}{Remarks}
\newtheorem{problems}{Problems}
\newtheorem{algorithms}{Algorithms}
\newtheorem{questions}{Questions}
\newtheorem{recalls}{Recalls}
\newtheorem{exercises}{Exercises}
\newtheorem{notes}{Notes}
\newtheorem{claims}{Claims}
\newtheorem{facts}{Facts}

\theoremstyle{definition}\newtheorem*{thm*}{Theorem}
\newtheorem*{lmm*}{Lemma}
\newtheorem*{dfn*}{Definition}
\newtheorem*{crl*}{Corollary}
\newtheorem*{prp*}{Proposition}
\newtheorem*{example*}{Example}
\newtheorem*{examples*}{Examples}
\newtheorem*{rmk*}{Remark}
\newtheorem*{rmks*}{Remarks}
\newtheorem*{prb*}{Problem}
\newtheorem*{alg*}{Algorithm}
\newtheorem*{recall*}{Recall}
\newtheorem*{recalls*}{Recalls}
\newtheorem*{exercise*}{Exercise}
\newtheorem*{exercises*}{Exercises}
\newtheorem*{note*}{Note}
\newtheorem*{notes*}{Notes}
\newtheorem*{caution*}{Caution}
\newtheorem*{claim*}{Claim}
\newtheorem*{fact*}{Fact}
\title[\textcolor{brown}{On Lipschitz Retraction}]{\textcolor{brown}{On Lipschitz Retraction of Finite Subsets of Normed Spaces}}

\author[\textcolor{gray}{E. Akofor}]{\textcolor{darkgray}{\\~\\ {\large Earnest Akofor \\~\\~\\ {\it Department of Mathematics \\~\\ Syracuse University}} \\~\\~\\ e-mail: eakofor@syr.edu}\\~\\}


\begin{document}

\begingroup
\let\MakeUppercase\relax 
\maketitle
\endgroup
\tableofcontents

\begin{abstract}
~\\~
\textcolor{darkgray}{If $X$ is a metric space, then its finite subset spaces $X(n)$ form a nested sequence under natural isometric embeddings $X = X(1)\subset X(2) \subset \cdots$. It was previously established, by Kovalev when $X$ is a Hilbert space and, by Ba\v{c}\'{a}k and Kovalev when $X$ is a CAT(0) space, that this sequence admits Lipschitz retractions $X(n)\rightarrow X(n-1)$ for all $n\geq 2$. We prove that when $X$ is a normed space, the above sequence admits Lipschitz retractions $X(n)\rightarrow X$, $X(n)\rightarrow X(2)$, as well as concrete retractions $X(n)\rightarrow X(n-1)$ that are Lipschitz if $n=2,3$ and H{\"o}lder-continuous on bounded sets if $n>3$. We also prove that if $X$ is a geodesic metric space, then each $X(n)$ is a $2$-quasiconvex metric space. These results are relevant to certain questions in the aforementioned previous work which asked whether Lipschitz retractions $X(n)\rightarrow X(n-1)$, $n\geq 2$, exist for $X$ in more general classes of Banach spaces.}
\end{abstract}
\section{Introduction}\label{Intro}

Given a metric space $X$, let $X(n)$ denote the collection of nonempty finite subsets of $X$ of cardinality at most $n$, viewed as a metric space with respect to the Hausdorff distance $d_H$. Then we have a nested sequence under natural isometric embeddings $X = X(1)\subset X(2) \subset \cdots$. If $X,Z$ are metric spaces, a map $f:X\ra Z$ is \emph{Lipschitz} if there is a number $\ld\geq 0$, (the least of) which we denote by $\Lip(f)$, such that $d(f(x),f(y))\leq \ld d(x,y)$ for all $x,y\in X$. If $Z\subseteq X$, a map $r:X\ra Z$ is a \emph{retraction} if its restriction to $Z$ is the identity map, i.e., $r|_Z=id_Z:Z\ra Z$, $z\mapsto z$. A map $X\ra Z$ is a \emph{Lipschitz retraction} if it is both a Lipschitz map and a retraction.

As pointed out by L. V. Kovalev in \cite{kovalev2016}, an example from \cite{mostovoy2004} shows that for a generic metric space $X$, there need not exist continuous retractions $r:X(n)\ra X(n-1)$ for every $n\geq 2$. That is, given $n\geq 2$, a continuous map of metric spaces $f:X(n-1)\ra Z$ need not always extend to a continuous map $F:X(n)\ra Z$. It was then proved in \cite{kovalev2016} that for a Hilbert space $\H$, however, Lipschitz retractions $r:\H(n)\ra \H(n-1)$ exist for all $n\geq 2$. This result was generalized by M. Ba\v{c}\'{a}k and L. V. Kovalev in \cite{bacac-kovalev2016} to the case where $X$ is a $\CAT(0)$ space.

In this paper, we prove that when $X$ is a normed space, the aforementioned sequence admits Lipschitz retractions $X(n)\ra X$, $X(n)\ra X(2)$, as well as concrete retractions $X(n)\ra X(n-1)$ that are Lipschitz if $n=2,3$ and H{\"o}lder-continuous on bounded sets if $n>3$. We also prove that if $X$ is a geodesic metric space, then each $X(n)$ is a $2$-quasiconvex metric space.

Question 3.4 of \cite{kovalev2016} asked whether Lipschitz retractions $X(n)\ra X(n-1)$ exist for all $n\geq2$ when $X$ is a Banach space. A related question in Remark 3.4 of \cite{bacac-kovalev2016} similarly asked whether strictly convex or uniformly convex Banach spaces admit such Lipschitz retractions. The results of this paper provide tools of investigation towards answering the above questions.

We begin in Section \ref{FinSub} by introducing our notation for finite subset spaces $X(n)$ of a metric space $X$. This involves Hausdorff distance, our definition of finite subset spaces of a metric space, and a few results we will use in subsequent sections.

In Section \ref{LipRetV3}, Theorem \ref{LrThmV3}, we show that if $X$ is a normed space, then there exists a Lipschitz retraction $X(3)\ra X(2)$ with Lipschitz constant $731$. In Section \ref{LipRetVn12}, Theorem \ref{LrThmVn2}, we use a result of P. Shvartsman in \cite[Theorem~1.2]{shvartsman2004} to construct Lipschitz retractions $X(n)\ra X$ and $X(n)\ra X(2)$. In Section \ref{Retractions}, Theorem \ref{MainThm}, using a refinement of the differentiable-path mapping technique applied in \cite{kovalev2016} when $X$ is a Hilbert space, we show that for any normed space $X$, there exist retractions $r:X(n)\ra X(n-1)$, for all $n\geq 2$, that are H{\"o}lder-continuous on bounded sets.

In Section \ref{QuasiConv}, we derive a key property of finite subset spaces $X(n)$ of a geodesic metric space $X$, namely, quasiconvexity. We show in Theorem \ref{QConvThm} that for any geodesic metric space $X$, the space $X(n)$ is $2$-quasiconvex for all $n\geq 1$, $X(2)$ is a geodesic metric space, and for $n\geq 3$ the quasiconvexity constant $2$ for $X(n)$ is the smallest possible.

\section{Finite Subset Spaces of Metric Spaces}\label{FinSub}
Throughout, we denote the collection of all subsets of a topological space $X$ by $\P(X)$, and nonempty subsets of $X$ by $\P^\ast(X)$. The \emph{closure} of a set $A\subseteq X$ will be denoted by $\ol{A}$. If $X$ is a metric space, $\B(X)$ and $\B^\ast(X)$ will denote bounded subsets and nonempty bounded subsets respectively. Furthermore, if $A,B\subseteq X$ and $\vep>0$, we will denote the \emph{distance} between $A,B$ by $\dist(A,B):=\inf\{d(a,b):a\in A,b\in B\}$, the \emph{open $\vep$-neighborhood} of $A$ by $N_\vep(A)=\{x\in X:\dist(x,A)<\vep\}$, and the \emph{closed $\vep$-neighborhood} of $A$ by $\ol{N}_\vep(A)=\{x\in X:\dist(x,A)\leq\vep\}$. The \emph{cardinality} of a set $A$ will be denoted by $|A|$.

\begin{dfn}[\textcolor{OliveGreen}{Hausdorff distance}] Let $X$ be a metric space. The Hausdorff distance on $\P^\ast(X)$ is the map $d_H:\P^\ast(X)\times\P^\ast(X)\ra\Real$ given by
\begin{align*}
d_H(A,B)&:=\max\left\{\sup_{a\in A}\inf_{b\in B}d(a,b)~,~\sup_{b\in B}\inf_{a\in A}d(a,b)\right\}.
\end{align*}
\end{dfn}
An alternative expression for $d_H$ is as follows.
\begin{rmk}\label{HausDistRmk}
Let $X$ be a metric space and $A,B\subseteq X$. Define $\ol{A}_\vep:=\ol{N}_\vep(A)$ and $\ol{B}_\vep:=\ol{N}_\vep(B)$, i.e., the closed $\vep$-neighborhoods of $A$ and $B$ in $X$. Then
\begin{equation*}
d_H(A,B)=\inf\left\{\vep:A\subseteq\ol{B}_\vep,~B\subseteq\ol{A}_\vep\right\}=\inf\left\{\vep:A\cup B\subseteq\ol{A}_\vep\cap\ol{B}_\vep\right\}.
\end{equation*}
In particular, if $A,B$ are compact, then $d_H(A,B)$ is ``achieved'' in the sense that with $\rho:=d_H(A,B)$, we have
\begin{equation*}
A\cup B\subseteq\ol{A}_\rho\cap\ol{B}_\rho.
\end{equation*}
\end{rmk}
A discussion of the basic properties of $d_H$ stated above can be found in \cite{BBI}, Section 7.3.1 (page 252).

\begin{dfn}[\textcolor{OliveGreen}{Symmetric product/power of a metric space}]\label{SymmProd}
Let $X$ be a metric space. For $n\geq 1$, the $n$th symmetric product (or the $n$th finite subset space) of $X$ is the set
\begin{equation*}
X(n):=\big\{A\subseteq X:|A|\leq n\big\}=\Big\{\{x_1,...,x_n\}:(x_1,...,x_n)\in X^n\Big\},
\end{equation*}
viewed as a metric space with respect to the Hausdorff distance
\begin{equation*}
d_H\big(\{x_1,...,x_n\},\{x_1',...,x_n'\}\big):=\max\left\{\max_{1\leq i\leq n}\min_{1\leq j\leq n}d(x_i,x_j'),\max_{1\leq i\leq n}\min_{1\leq j\leq n}d(x_j,x_i')\right\}.
\end{equation*}
\end{dfn}

\begin{dfn}[\textcolor{OliveGreen}{Minimum separation}]\label{MinSep}
Let $X$ be a metric space. The minimum separation is the function $\delta:X(n)\ra[0,+\infty)$ given by $\delta(x)=\min\limits_{i\neq j}d(x_i,x_j)$.
\end{dfn}
Note that for any $x\in X(n)$, if $|x|<n$ then $\delta(x)=0$, i.e., $\delta|_{X(n-1)}=0$.

\begin{lmm}\label{HausDistBound}
Let $X$ be a metric space. For all $x,y\in X(n)$, if $\delta(x)>2d_H(x,y)$ or $\delta(y)>2d_H(x,y)$, then there exists an enumeration $x_i,y_i$ of elements of $x,y$ such that
\begin{equation}
\label{KOVeq8}d(x_i,y_i)\leq d_H(x,y),~~~~~~\txt{for all}~~i=1,...,n.
\end{equation}
\end{lmm}
\begin{proof}
By symmetry it suffices to assume $\delta(x)>2d_H(x,y)$. Let $\rho:=d_H(x,y)$. Then by Remark \ref{HausDistRmk}, $x\subseteq \ol{y}_\rho=\bigcup_{i=1}^n\ol{N}_\rho(y_i)$ and $y\subseteq \ol{x}_\rho=\bigcup_{i=1}^n\ol{N}_\rho(x_i)$. Moreover, the balls $\ol{N}_\rho(x_i)$ are disjoint because $\delta(x)>2\rho$. Thus, each ball $\ol{N}_\rho(x_i)$ contains exactly one point of $y$, i.e., for each $x_i\in x$, there exists a \emph{unique} $y_j\in y$ such that $d(x_i,y_j)\leq\rho$. This means we can \emph{relabel} the points $x_i,y_i$ so that (\ref{KOVeq8}) holds.
\end{proof}

\begin{lmm}\label{DeltaCont}
Let $X$ be a metric space. The minimum separation $\delta$ in Definition \ref{MinSep} is $2$-Lipschitz.
\end{lmm}
\begin{proof}
If $\delta(x)\leq 2d_H(x,y)$ and $\delta(y)\leq 2d_H(x,y)$, it is clear that
\begin{equation*}
|\delta(x)-\delta(y)|\leq\max\big(\delta(x),\delta(y)\big)\leq 2d_H(x,y).
\end{equation*}
So, assume $\delta(x)>2d_H(x,y)$ or $\delta(y)>2d_H(x,y)$. Then using Lemma \ref{HausDistBound},
\begin{equation*}
|\delta(x)-\delta(y)|\leq\max_{i\neq j}\Big|d(x_i,x_j)-d(y_i,y_j)\Big|\leq \max_{i\neq j}\Big(d(x_i,y_i)+d(x_j,y_j)\Big)\sr{(\ref{KOVeq8})}{\leq} 2d_H(x,y).\qedhere
\end{equation*}
\end{proof}

\begin{lmm}[\textcolor{OliveGreen}{Diameter is 2-Lipschitz}]\label{DiamCont}
If $X$ is a metric space, the map $\diam:\B^\ast(X)\ra \Real$ given by $\diam(A):=\sup\{d(a,a'):a,a'\in A\}$ is $2$-Lipschitz with respect to the Hausdorff distance $d_H$.
\end{lmm}
\begin{proof}
Fix $\vep>0$ and bounded sets $A,B\in \B^\ast(X)$. Then there exist $a,a'\in A$ and $b,b'\in B$ such that $\diam(A)\leq d(a,a')+\vep$, $d(a,b)\leq \dist(a,B)+\vep\leq d_H(A,B)+\vep$, and $d(a',b')\leq \dist(a',B)+\vep\leq d_H(A,B)+\vep$. These three inequalities (together with the triangle inequality) in turn imply
\begin{equation*}
\diam(A)\leq d(a,b)+d(b,b')+d(a',b')+\vep\leq 2d_H(A,B)+\diam(B)+3\vep.
\end{equation*}
Hence, $|\diam(A)-\diam(B)|\leq 2 d_H(A,B)+3\vep$.
\end{proof}

\section{Concrete Lipschitz Retractions $X(2)\ra X$ and $X(3)\ra X(2)$}\label{LipRetV3}
The Lipschitz retractions in this section, unlike those in Section \ref{LipRetVn12}, have concrete Lipschitz constants.

\begin{dfn}[\textcolor{OliveGreen}{Addition and scalar multiplication of sets}]\label{SumSclMult} Let $X$ be a vector space. If $A,B\subseteq X$ and $\ld$ is a scalar, we write $A+B:=\{a+b:a\in A,b\in B\}$ and $\ld A:=\{\ld a:a\in A\}$.
\end{dfn}

\begin{dfn}[\textcolor{OliveGreen}{Scale-invariant map, Translation-invariant map, Affine map}]
Let $X$ be a vector space and $f:\E\subseteq\P^\ast(X)\ra\P^\ast(X)$ a map. We say $f$ is scale-invariant (or commutes with scaling) if for any $t\in \Real$, we have $f(tA)=tf(A)$ for all $A\in\E$ such that $tA\in\E$. Similarly, $f$ is translation-invariant (or commutes with translations) if for any $v\in X$, we have $f(A+v)=f(A)+v$ for all $A\in\E$ such that $A+v\in\E$.

We say $f$ is affine if $f$ is both scale-invariant and translation-invariant, i.e., for any $t\in\Real$, $v\in X$, we have $f(tA+v)=tf(A)+v$ for all $A\in\E$ such that $tA+v\in\E$.
\end{dfn}

\begin{dfn}[\textcolor{OliveGreen}{Proximal map between points of $X(n)$}]
Let $X$ be a metric space and $x,y\in X(n)$. A map $p:x\ra y$ is proximal if $d(a,p(a))\leq d_H(x,y)$ for all $a\in x$.
\end{dfn}

\begin{dfn}[\textcolor{OliveGreen}{Normalized element in $X(n)$, Set of normalized elements}]
Let $X$ be a metric space and $x\in X(n)$. We say $x$ is normalized if $\diam(x)=1$. We will write $N\big(X(n)\big)$ for all normalized elements of $X(n)$.
\end{dfn}

\begin{dfn}[\textcolor{OliveGreen}{Central element in $X(n)$, Set of central elements, Set of normalized central elements}]
Let $X$ be a normed space and $x\in X(n)$. We say $x$ is central if $0\in x$. We will write $X_0(n)=\{x\in X(n):0\in x\}$ for all central elements of $X(n)$. Accordingly, we will write $N(X_0(n))$ for the set of normalized central elements of $X(n)$.
\end{dfn}

Note that every element $x\in X(n)$ can be written (not uniquely) as
\begin{equation*}
x=t x_0+v,~~~~\txt{for some}~~t\in[0,+\infty),~~x_0\in N(X_0(n)),~~v\in X.
\end{equation*}

\begin{lmm}[\textcolor{OliveGreen}{Homogeneous Lipschitz Extension}]\label{HomExtLmm}
Let $X$ be a normed space and $1\leq k\leq n-1$. Any translation-invariant Lipschitz map $R:N(X_0(n))\ra X(k)$ satisfying $R(x)\subseteq\Conv(x)$ and $R|_{N(X_0(n))\cap X(k)} = id$ can be extended to an affine Lipschitz retraction $r:X(n)\ra X(k)$ with $\Lip(r)=6\Lip(R)+5$.
\end{lmm}
\begin{proof}
Let $R:N(X_0(n))\ra X(k)$ be a Lipschitz map such that $R(x)\subseteq\Conv(x)$ and $R|_{N(X_0(n))\cap X(k)}$ = $id$. Define a map $r:X(n)\ra X(k)$ by $r(tx+v)=tR(x)+v$ for all $x\in N(X_0(n))$, $t\in[0,+\infty)$, and $v\in X$. Then $r$ is well defined because $R$ is translation-invariant. For any $x,y\in N(X_0(n))$ and $t,s\in[0,+\infty)$, since $0\in x$ and diameter is $2$-Lipschitz with respect to $d_H$,
\begin{equation*}
d_H(tx,sx)\leq |t-s|=|\diam(tx)-\diam(s y)|\leq 2d_H(tx,sy).
\end{equation*}
 Thus, using the triangle inequality, $0\in y$, and $R(y)\subseteq \Conv(y)$, we have
\begin{align*}
d_H(&r(tx),r(sy))\leq d_H(r(tx),r(ty))+d_H(r(ty),r(sy))=d_H(tR(x),tR(y))+d_H(tR(y),sR(y))\\
    &\leq\Lip(R)d_H(tx,ty)+\diam(R(y))|t-s|\leq \Lip(R)\Big[d_H(tx,sy)+d_H(sy,ty)\Big]+|t-s|\\
    &\leq\Lip(R)d_H(tx,sy)+\Big(\Lip(R)+1\Big)|t-s|\leq \Big(3\Lip(R)+2\Big)d_H(tx,sy),
\end{align*}
which shows $r$ is Lipschitz on $X_0(n)$. Given $x,y\in X(n)$, let $u\in x$, $v\in y$ such that $\|u-v\|\leq d_H(x,y)$. Then
\begin{align*}
d_H\big(r(x),r(y)\big)&=d_H\big(r(x-u)+u,r(y-v)+v\big)\leq d_H\big(r(x-u),r(y-v)\big)+\|u-v\|\\
 &\leq\big(2\Lip(r|_{X_0(n)})+1\big)d_H(x,y).\qedhere
\end{align*}
\end{proof}

\begin{dfn}[\textcolor{OliveGreen}{Thin sets in $X(3)$}]
Let $X$ be a normed space and let $x\in X(3)$. Then $x$ is called ``thin in $X(3)$'' if $x$ is normalized and $0\leq\delta(x)\leq{1\over 3}$. We will denote all thin sets in $X(3)$ by $\Thin\big(X(3)\big)$.
\end{dfn}

\begin{notation}
If $x=\{x_1,x_2,x_3\}\in \Thin\big(X(3)\big)$, we will assume without loss of generality that
\begin{equation*}
\delta(x)=d(x_1,x_2)\leq d(x_2,x_3)\leq d(x_1,x_3)=1,
\end{equation*}
where the triangle inequality implies $d(x_2,x_3)\geq 2/3$.
\end{notation}

\begin{dfn}[\textcolor{OliveGreen}{Vertex map of thin elements of  $X(3)$}]\label{ExtrMinDfn}
This is the map $V:\Thin\big(X(3)\big)\ra X$ given by $V(\{x_1,x_2,x_3\}):=x_3$.
\end{dfn}

\begin{dfn}[\textcolor{OliveGreen}{Average map}]
If $X$ is a normed space, the average $\Avg:X(n)\ra X$ is given by $\Avg(x)~:=~{1\over |x|}\sum_{a\in x}a$. In particular, if $x=\{x_1,...,x_n\}\in X(n)\backslash X(n-1)$, then we can write $\Avg(x)={1\over n}\sum_{i=1}^n x_i$.
\end{dfn}

Throughout the rest of this section, we will assume $X$ is a normed space.

\begin{lmm}[\textcolor{OliveGreen}{Continuity of the Average map}]\label{ContAvgLmm}
Let $X$ be a normed space and $x,y\in X(n)\backslash X(n-1)$.
\bit
\item[(i)] If a proximal bijection $x\ra y$ exists, then ~~$d_H(\Avg(x),\Avg(y))\leq d_H(x,y)$.
\item[(ii)] If $\max\Big(\delta(x),\delta(y)\Big)\leq 2d_H(x,y)$, and $\al\diam(x)\leq\delta(x)$ or $\al\diam(y)\leq\delta(y)$ for constant $\al>0$, then
\begin{equation*}
\textstyle d_H(\Avg(x),\Avg(y))\leq \left(1+{2\over\al}\right)d_H(x,y).
\end{equation*}
\eit
\end{lmm}
\begin{proof}
(i) Let $x\ra y$, $x_i\mapsto y(i)$ be a proximal bijection. Then
\begin{align*}
\textstyle d_H(\Avg(x),\Avg(y))=\left\|{\sum x_i\over n}-{\sum y_i\over n}\right\|=\left\|{\sum x_i\over n}-{\sum y(i)\over n}\right\|\leq d_H(x,y).
\end{align*}

(ii) It suffices to assume $\al\diam(y)\leq\delta(y)$. Consider a proximal map $x\ra y$, $x_i\mapsto y(i)$. Then
\begin{align*}
d_H(&\Avg(x),\Avg(y))\textstyle=\left\|{\sum x_i\over n}-{\sum y_i\over n}\right\|\leq {\sum\|x_i-y(i)\|+\sum_i\|y(i)-y_i\|\over n}\\
  &\textstyle\leq d_H(x,y)+\diam(y)\leq d_H(x,y)+{\delta(y)\over\al}\leq\left(1+{2\over \al}\right)d_H(x,y).\qedhere
\end{align*}
\end{proof}

\begin{prp}\label{LrPrpV2}
Let $X$ be a normed space. There exists a $1$-Lipschitz retraction $X(2)\ra X$.
\end{prp}
\begin{proof}
Define $r:X(2)\ra X$ by $r(x)=\Avg(x)$. Then $r|_X=id$. If $x,y\in X(n)$, consider the following.
\bit
\item[(i)] $x=\{x_1\}\in X$ and $y=\{y_1,y_2\}\in X(2)\backslash X$:~ In this case,
\begin{equation*}
\|r(x)-r(y)\|\leq \max(\|x_1-y_1\|,\|x_1-y_2\|)=d_H(x,y).
\end{equation*}
\item[(ii)] $x,y\in X(2)\backslash X$: A proximal bijection $x\ra y$ exists, and so by Lemma \ref{ContAvgLmm}(i), $d_H(r(x),r(y))\leq d_H(x,y)$.
\eit
Hence, $r$ is a $1$-Lipschitz retraction.
\end{proof}

\begin{lmm}\label{LipCenLmm1}
If $X$ is a normed space, the following map is $3$-Lipschitz.
\begin{equation*}
\textstyle f:\Thin\big(X(3)\big)\ra X(2),~~~~f(x):=\Big\{\Avg\big(x\backslash V(x)\big),V(x)\Big\}=\left\{{x_1+x_2\over 2},x_3\right\}.
\end{equation*}
\end{lmm}
\begin{proof}
Let $x,y\in X(3)$ be thin sets. Observe that $d_H(f(x),x)\leq{1\over 2}\delta(x)$, and so
\begin{equation*}
\textstyle d_H(f(x),f(y))\leq {1\over 2}\delta(x)+{1\over 2}\delta(y)+d_H(x,y).
\end{equation*}
Thus, if $\delta(x)\leq 2d_H(x,y)$ and $\delta(y)\leq 2d_H(x,y)$, then $d_H(f(x),f(y))\leq 3d_H(x,y)$. So, assume
\begin{equation*}
\textstyle d_H(x,y)<{1\over 2}\delta(x)~~~~\txt{or}~~~~ d_H(x,y)<{1\over 2}\delta(y),~~~~~~~~\left(\Ra~~d_H(x,y)<{1\over 2}{1\over 3}={1\over 6}\right).
\end{equation*}
Then by Lemma \ref{HausDistBound}, we have a proximal bijection $x_i\mapsto y(i)$ such that $\|x_i-y(i)\|\leq d_H(x,y)<1/6$ for all $i$. This bijection satisfies $\{x_1,x_2\}\mapsto\{x(1),x(2)\}=\{y_1,y_2\}$, i.e., $x_3\mapsto y(3)=y_3$, since
\begin{equation*}
\|y(1)-y(2)\|\leq\|x_1-y(1)\|+\|x_1-x_2\|+\|x_2-y(2)\|<{2\over 3}\leq\|y_2-y_3\|.
\end{equation*}
Hence,
\begin{equation*}
\textstyle d_H(f(x),f(y))\leq \max\left(\left\|{x_1+x_2\over 2}-{y_1+y_2\over 2}\right\|,\|x_3-y_3\|\right)\leq d_H(x,y).\qedhere
\end{equation*}
\end{proof}

\begin{dfn}[\textcolor{OliveGreen}{Lipschitz partition of unity}]\label{LipPoU}
Consider the maps $\vphi_1,\vphi_2:\Real\ra\Real$ given by
\begin{equation*}
\vphi_1(t)=
\left\{
  \begin{array}{ll}
    1, & t\leq {1\over 5} \\
    -20t+5, & t\in \left[{1\over 5},{1\over 4}\right]\\
    0, & t\geq {1\over 4}
  \end{array}
\right\},~~~~
\vphi_2(t)=
\left\{
  \begin{array}{ll}
    0, & t\leq {1\over 5}\\
    20t-4, & t\in\left[{1\over 5},{1\over 4}\right]\\
    1, & t\geq {1\over 4}
  \end{array}
\right\}.
\end{equation*}
\end{dfn}
The above maps form a $20$-Lipschitz partition of unity.

\begin{lmm}[\textcolor{OliveGreen}{Gluing with strips}]\label{LSGluLmm}
Let $X$ be a metric space and $\vphi:X\ra\Real$ a Lipschitz function. Consider a finite collection of intervals $\{I_k=(a_k,b_k):k=1,...,N\}$ such that $a_k<a_{k+1}<b_k<b_{k+1}$ for all $k=1,...,N-1$ and $\Real=\bigcup_{k=1}^N(a_k,b_k)$.

If a map $g:E\subseteq X\ra X$ satisfies $\sup_{x\in E}d(x,g(x))\leq D<\infty$ and is Lipschitz on each of the sets $E_k=\vphi^{-1}(a_k,b_k)=\{x\in E:a_k<\vphi(x)<b_k\}$, then $g$ is Lipschitz, and
\begin{equation*}
\textstyle\Lip(g)=\max\left\{\max\limits_{1\leq k\leq N}\Lip(g|_{E_k}),{\left(1+{2D\Lip(\vphi)\over\vep}\right)}\right\},
\end{equation*}
where ~$\vep:=\min\limits_k\diam(I_k\cap I_{k+1})=\min\limits_{1\leq k\leq N-1}|a_{k+1}-b_k|$.
\end{lmm}
\begin{proof}
Pick any $x,y\in E$. If $|\vphi(x)-\vphi(y)|<\vep$, then $x,y\in E_k=\vphi^{-1}(a_k,b_k)$ for some $k$, and so $d(g(x),g(y))\leq \Lip(g|_{E_k})d(x,y)$. On the other hand, if $|\vphi(x)-\vphi(y)|\geq\vep$, then by the triangle inequality, we get
\begin{equation*}
\textstyle d(g(x),g(y))\leq \left(1+{2D\Lip(\vphi)\over\vep}\right)d(x,y).\qedhere
\end{equation*}
\end{proof}

\begin{dfn}[\textcolor{OliveGreen}{Interpolation map of $X(3)$}]\label{IntMapDfnV3}
This is the map $R:N(X_0(3))\ra X(2)$ given by
\begin{equation*}
R(x)=\vphi_1\big(\delta(x)\big)R_1(x)+\vphi_2\big(\delta(x)\big)R_2(x),
\end{equation*}
where $\vphi_1,\vphi_2$ are as in Definition \ref{LipPoU},
\begin{equation*}
\textstyle R_1(x):=f(x)=\left\{{x_1+x_2\over 2},x_3\right\},~~~~R_2(x):=\Avg(x)={x_1+x_2+x_3\over 3},
\end{equation*}
and, we add and multiply sets by scalars as in Definition \ref{SumSclMult}.
\end{dfn}

\begin{lmm}\label{LipIntLmm}
The interpolation map $R:N(X_0(3))\ra X(2)$ is $121$-Lipschitz.
\end{lmm}
\begin{proof}
By Lemma \ref{LSGluLmm} with $\vphi=\delta$, $g=R$, $\{(a_k,b_k)\}=\{(-\infty,1/5),(1/6,1/3),(1/4,+\infty)\}$, it suffices to show that $R$ is Lipschitz on each of the following sets (within $N(X_0(3))$).
\begin{equation*}
\textstyle A=\left\{\delta\leq {1\over 5}\right\},~~~~C=\left\{{1\over 6}\leq\delta\leq{1\over 3}\right\} ,~~~~B=\left\{\delta\geq{1\over 4}\right\}.
\end{equation*}

\ul{$R$ is $3$-Lipschitz on $A$}: This follows from Lemma \ref{LipCenLmm1}.

\ul{$R$ is $9$-Lipschitz on $B$}: If $x,y\in B$, then $1\leq 4\delta(x),4\delta(y)$. Consider cases as follows. If $\delta(x)>2d_H(x,y)$ or $\delta(y)>2d_H(x,y)$, we have a proximal bijection $x\ra y$, and so by Lemma \ref{ContAvgLmm}(i),
\begin{equation*}
d_H(R(x),R(y))\leq d_H(x,y).
\end{equation*}
On the other hand, if $\delta(x),\delta(y)\leq 2d_H(x,y)$, then by Lemma \ref{ContAvgLmm}(ii),
\begin{equation*}
\textstyle d_H(R(x),R(y))\leq\left(1+{2\over 1/4}\right)d_H(x,y)=9d_H(x,y).
\end{equation*}

\ul{$R$ is $44$-Lipschitz on $C$}: This follows from the fact that $\vphi_i(\delta(x))$, $R_i(x)$, $i=1,2$ are bounded Lipschitz maps. With $\delta=\delta(x)$ and $\delta'=\delta(y)$, we have
{\small\begin{align*}
d_H(R(x)&,R(y))=d_H\Big(\vphi_1(\delta)R_1(x)+\vphi_2(\delta)R_2(x)~,~\vphi_1(\delta')R_1(y)+\vphi_2(\delta')R_2(y)\Big)\\
    &\leq d_H\Big(\vphi_1(\delta)R_1(x)~,~\vphi_1(\delta')R_1(y)\Big)+\Big\|\vphi_2(\delta)R_2(x)-\vphi_2(\delta')R_2(y)\Big\|\\
    &\leq \big[(20+3)+(20+1)\big]d_H(x,y)=44d_H(x,y).
\end{align*}}
The Lipschitz constant of $R$ can be calculated from Lemma \ref{LSGluLmm} as
\begin{equation*}
\textstyle \Lip(R)=\max\left\{44,1+{2\times 1\times 2\over 1/5-1/6}\right\}=121.\qedhere
\end{equation*}
\end{proof}

\begin{thm}\label{LrThmV3}
Let $X$ be a normed space. There exist Lipschitz retractions $X(n)\ra X(n-1)$ for $n=2,3$. From Lemmas \ref{HomExtLmm} and \ref{LipIntLmm}, the retraction $X(3)\ra X(2)$ has Lipschitz constant $6(121)+5=731$.
\end{thm}
\begin{proof}
The case of $n=2$ is given by Proposition \ref{LrPrpV2}. So, let $n=3$. Then by Lemma \ref{HomExtLmm}, it is enough to show that the interpolation map $R:N(X_0(3))\ra X(2)$ is Lipschitz, which follows from Lemma \ref{LipIntLmm}.
\end{proof}

\section{Lipschitz Retractions $X(n)\ra X$ and $X(n)\ra X(2)$}\label{LipRetVn12}

If not stated, we will assume $X$ is a normed space. Let $\K_n(X)$:=$\{$convex compact subsets of $X$ of dimension $\leq n$$\}$. By Theorem 1.2 of \cite{shvartsman2004}, there exists an affine Lipschitz selector $S:\K_n(X)\ra X,~A\mapsto S(A)\in A$.
\begin{prp}\label{SvartLRVn1}
If $X$ is a normed space, there exist affine Lipschitz retractions $X(n)\ra X$ for all $n\geq 1$.
\end{prp}
\begin{proof}
Consider the map $s=S\circ\Conv:X(n)\sr{\Conv}{\ral}\K_n(X)\sr{S}{\ral}X$, where the convex hull operation is
\begin{equation*}
\textstyle \Conv(x)=\left\{\sum_{i=1}^n\al_ix_i:~\sum_{i=1}^n\al_i=1,~\al_i\geq0\right\}.
\end{equation*}
Given $\sum\al_ix_i\in \Conv(x)$, it follows from the definition of $d_H(x,y)$ that for each $x_i\in x$, there exists $y(i)\in y$ such that $\|x_i-y(i)\|\leq d_H(x,y)$. Since $\sum\al_iy(i)\in\Conv(y)$ and $\left\|\sum\al_i x_i-\sum\al_iy(i)\right\|\leq d_H(x,y)$, it follows by symmetry in the definition of Hausdorff distance that $d_H(\Conv(x),\Conv(y))\leq d_H(x,y)$. Hence,
\begin{align*}
d_H(s(x),s(y))&=\|s(x)-s(y)\|=\|S\circ \Conv(x)-S\circ \Conv(y)\|\\
 &\leq\Lip(S)d_H(x,y).\qedhere
\end{align*}
\end{proof}

\begin{dfn}[\textcolor{OliveGreen}{Two-cluster decomposition of an element of $X(n)$}]
Let $x\in X(n)$ and consider numbers $\al,\beta>0$. A decomposition $x=x'\cup x''$ is an $(\al,\beta)$-decomposition if $x',x''$ are nonempty, $\diam(x')\leq\al$, $\diam(x'')\leq\al$, and $\dist(x',x'')\geq\beta$.
\end{dfn}

\begin{rmk}\label{UniqDecRmk}
Observe that an $(\al,\beta)$-decomposition $x=x'\cup x''$ is unique (up to permutation of the clusters $x',x''$) if $\al<\beta$. This is because if $x=\wt{x}'\cup\wt{x}''$ is any $(\al,\beta)$-decomposition, then neither $\wt{x}'$ nor $\wt{x}''$ can intersect both $x'$ and $x''$.

Moreover,  if $\alpha\le \alpha'<\beta' \le \beta$, then the $(\alpha', \beta')$-decomposition is the same as the $(\alpha,\beta)$-decomposition. This is because the $(\alpha,\beta)$-decomposition is also an $(\alpha', \beta')$-decomposition, which is unique. In particular, if $x$ has an $(\al,\beta)$-decomposition with $\al<\beta$, then for any number $0<c<{\beta-\al\over 2}$, any $(\al+c,\beta-c)$-decomposition is unique and equals the $(\al,\beta)$-decomposition.
\end{rmk}

Fix a number $\tau>6$.

\begin{dfn}[\textcolor{OliveGreen}{2nd order thin sets in $X(n)$, Collection of thin sets}]
Let $X$ be a metric space and $x\in X(n)$. We say $x$ is a thin set of order $2$ (or $2$-thin set) if $x$ is normalized  (i.e., $\diam(x)=1$) and ~{\small $\dist_H(x,X(2)):=\inf\limits_{z\in X(2)}d_H(x,z)<{1\over\tau}$}. We will denote the collection of all 2-thin sets in $X(n)$ by $\Thin_2\big(X(n)\big)$.
\end{dfn}
\begin{lmm}[\textcolor{OliveGreen}{Cluster decomposition of a $2$-thin set}]
Let $X$ be a normed space and $x\in\Thin_2(X(n))$. Then $x$ admits a unique $\left({2\over\tau},1-{4\over\tau}\right)$-decomposition.
\end{lmm}
\begin{proof}
(i) \ul{Existence}: Since $\dist_H(x,X(2))<{1\over\tau}$, there exists $\{a,b\}\in X(2)$ such that $d_H(x,\{a,b\})<{1\over\tau}$. Since diameter is 2-Lipschitz with respect to $d_H$, we have $|\diam(x)-\diam(\{a,b\})|\leq 2d_H(x,\{a,b\})$. Thus,
\begin{equation*}
\textstyle \|a-b\|>1-{2\over\tau}.
\end{equation*}
Observe that for any $u\in x$, we have either $\|u-a\|<{1\over\tau}$ or $\|u-b\|<{1\over\tau}$ but not both: Otherwise, if  $\|u-a\|<{1\over\tau}$ and $\|u-b\|<{1\over\tau}$ then the triangle inequality gives $1-{2\over\tau}<\|a-b\|<{2\over\tau}$, which implies $\tau<4$ (a contradiction since $\tau>6$ by assumption). Let
\begin{align*}
x'&:=\{u\in x:\|u-a\|<1/\tau\}=x\cap N_{1/\tau}(a),\\
x''&:=\{u\in x:\|u-b\|<1/\tau\}=x\cap N_{1/\tau}(b).
\end{align*}
Note that for any $u\in x'$, $v\in x''$, we have $\big|\|u-v\|-\|a-b\|\big|<{2\over\tau}$, and so $\|u-v\|>1-{4\over\tau}$. Hence,
\begin{equation}
\label{2ThinDefEq}\textstyle \diam(x')<{2\over\tau},~~~~\diam(x'')<{2\over\tau},~~~~\dist(x',x'')\geq 1-{4\over\tau},
\end{equation}
where $x',x''$ are nonempty because $\diam(x)=1>{2\over\tau}$.

(ii) \ul{Uniqueness}: This follows from Remark \ref{UniqDecRmk} and the fact that $\tau>6$.
\end{proof}

\begin{dfn}[\textcolor{OliveGreen}{Skeleton map}]\label{SkelMapDfn}
This is the map $J:N\big(X(n)\big)\ra X(2)$ with
\begin{equation}
\label{SkelMapEq}J(x):=\left\{
       \begin{array}{ll}
         R_1(x):=\Big\{s(x')~,~s(x'')\Big\}, & x\in \Thin_2\big(X(n)\big) \\
       R_2(x):=\Big\{s(x)\Big\}, & x\in N\big(X(n)\big)\backslash \Thin_2\big(X(n)\big)\\
       \end{array}
     \right\},
\end{equation}
where $x\in\Thin_2\big(X(n)\big)$ decomposes as $x=x'\cup x''$, and $s:X(n)\ra X$ are the affine Lipschitz retractions from Proposition \ref{SvartLRVn1}.
\end{dfn}

\begin{lmm}\label{FirstLipLmm}
The map $R_1:\Thin_2\big(X(n)\big)\ra X(2)$ is Lipschitz.
\end{lmm}
\begin{proof}
If $x\in\Thin_2(X(n))$, let $x=x'\cup x''$ be the $(\al,\beta)=(2/\tau,1-4/\tau)$-decomposition of $x$. Pick a number $0<\rho<{\beta-\al\over 4}={1\over 4}\left(1-{6\over\tau}\right)$. Let $y\in\Thin_2(X(n))$ such that $d_H(x,y)\leq\rho$, and define
\begin{align*}
y'&=\{u\in y:\dist(u,x')\leq\rho\}=y\cap N_\rho(x'),\\
y''&=\{u\in y:\dist(u,x'')\leq\rho\}=y\cap N_\rho(x'').
\end{align*}
Observe that $y=y'\cup y''$ (by the definition of Hausdorff distance), $\diam(y')<\al+2\rho$, $\diam(y'')<\al+2\rho$, and $\dist(y',y'')>\beta-2\rho$. Thus, $y',y''$ give a unique $(\al+2\rho,\beta-2\rho)$-decomposition since $\al+2\rho<\beta-2\rho$ (where $\al+2\rho<{1\over 2}-{1\over\tau}<1$, and so $y',y''$ are nonempty). By Remark \ref{UniqDecRmk}, this $(\al+2\rho,\al-2\rho)$-decomposition of $y$ is the same as the $(\al,\beta)$-decomposition of $y$.

By construction, $x'\cup y'$ and $x''\cup y''$ are \emph{distantly separated} in the sense that
\begin{equation*}
\dist(x'\cup y',x''\cup y'')>\beta-2\rho>\al+2\rho>\max\left(\diam(x'\cup y'),\diam(x''\cup y'')\right).
\end{equation*}
It follows by direct calculation that $d_H(x,y)=\max\left(d_H(x',y'),d_H(x'',y'')\right)$. Hence, we have
\begin{align*}
d_H(R_1(x),R_1(y))&=d_H\big(\{s(x'),s(x'')\},\{s(y'),s(y'')\}\big)\leq \Lip(s)\max\Big(d_H(x',y'),d_H(x'',y'')\Big)\\
  &\leq \Lip(s)d_H(x,y),~~~~\txt{if}~~~~d_H(x,y)\leq\rho.
\end{align*}

On the other hand, since $R_1(x)\subseteq \Conv(x)$, we also have
\begin{align*}
d_H(R_1(x)&,R_1(y))\leq d_H(x,y)+d_H(R_1(x),x)+d_H(R_1(y),y)\\
  &\leq d_H(x,y)+\diam(x)+\diam(y)=d_H(x,y)+2\\
  &\textstyle\leq \left(1+{2\over\rho}\right)d_H(x,y),~~~~\txt{if}~~~~d_H(x,y)\geq\rho.\qedhere
\end{align*}
\end{proof}

\begin{dfn}[\textcolor{OliveGreen}{Lipschitz partition of unity}]\label{LipPoUVn2}
Fix $\tau>0$. The functions $\vphi_1,\vphi_2:\Real\ra\Real$ given by
\begin{equation*}
\vphi_1(t)=
\left\{
  \begin{array}{ll}
    1, & t\leq{1\over 3\tau} \\
    -(6\tau)t+3, & t\in\left[{1\over 3\tau},{1\over 2\tau}\right]\\
    0, & t\geq {1\over 2\tau}
  \end{array}
\right\},~~~~
\vphi_2(t)=
\left\{
  \begin{array}{ll}
    0, & t\leq{1\over 3\tau} \\
    (6\tau)t-2, & t\in\left[{1\over 3\tau},{1\over 2\tau}\right]\\
    1, & t\geq {1\over 2\tau}
  \end{array}
\right\}
\end{equation*}
form a Lipschitz partition of unity.
\end{dfn}

\begin{dfn}[\textcolor{OliveGreen}{Interpolation map}]
This is the map $R:N(X_0(n))\ra X(2)$ given by
\begin{equation*}
R(x)=\vphi_1(\delta)R_1(x)+\vphi_2(\delta)R_2(x),~~~~\delta:=\dist_H\big(x,X(2)\big),
\end{equation*}
where $R_1,R_2$ are as in (\ref{SkelMapEq}), $\vphi_1,\vphi_2$ are as in Definition \ref{LipPoUVn2} and, we add and multiply sets by scalars as in Definition \ref{SumSclMult}.
\end{dfn}

\begin{lmm}\label{LipIntLmmVn2}
The interpolation map $R:N(X_0(n))\ra X(2)$ is Lipschitz.
\end{lmm}
\begin{proof}
By Lemma \ref{LSGluLmm} with $\vphi=\delta$, $g=R$, $\{(a_k,b_k)\}$ = $\{(-\infty,1/(3\tau))$, $(1/(4\tau),1/\tau)$, $(1/(2\tau),+\infty)\}$, it suffices to show that $R$ is Lipschitz on each of the following sets (within $N(X_0(n))$).
\begin{equation*}
\textstyle A=\left\{\delta\leq {1\over 3\tau}\right\},~~~~C=\left\{{1\over 4\tau}\leq\delta\leq {1\over\tau}\right\},~~~~B=\left\{\delta\geq{1\over 2\tau}\right\}.
\end{equation*}
$R$ is Lipschitz on $A$ by Lemma \ref{FirstLipLmm}, and Lipschitz on $B$ by the definition of $s$. On $C$, with $\delta=\dist_H\big(x,X(2)\big)$ and $\delta'=\dist_H\big(y,X(2)\big)$, we have
\begin{align*}
d_H(R(x)&,R(y))=d_H\Big(\vphi_1(\delta)R_1(x)+\vphi_2(\delta)R_2(x),\vphi_1(\delta')R_1(y)+\vphi_2(\delta')R_2(y)\Big)\\
    &\leq d_H\Big(\vphi_1(\delta)R_1(x),\vphi_1(\delta')R_1(y)\Big)+\Big\|\vphi_2(\delta)R_2(x)-\vphi_2(\delta')R_2(y)\Big\|.
\end{align*}
The result now follows because $\vphi_i(\delta)$, $R_i\big(x\big)$, $i=1,2$ are bounded Lipschitz maps.
\end{proof}

\begin{thm}\label{LrThmVn2}
Let $X$ be a normed space. There exist Lipschitz retractions $X(n)\ra X$ and $X(n)\ra X(2)$.
\end{thm}
\begin{proof}
The case of $X(n)\ra X$ is Proposition \ref{SvartLRVn1}. So, consider the case of $X(n)\ra X(2)$. By Lemma \ref{HomExtLmm}, it is enough to show that the interpolation map $R:N(X_0(n))\ra X(2)$ is Lipschitz, which follows from Lemma \ref{LipIntLmmVn2}.
\end{proof}

\section{Concrete H{\"o}lder Retractions $X(n)\ra X(n-1)$}\label{Retractions}
In this section, unless stated otherwise, $X$ is a real normed space. If $x\in X$ and $x^\ast\in X^\ast$, we will sometimes write the number $x^\ast(x)$ as $\langle x,x^\ast\rangle$ for convenience.

\begin{dfn}[\textcolor{OliveGreen}{Norming functional}] A linear functional $x^\ast\in X^\ast$ is a norming functional for $x_0\in X$ if $\|x^\ast\|=1$ and $x^\ast(x_0)=\|x_0\|$. If $x^\ast$ is a norming functional of $x_0$, we will also refer to $z^\ast:=\|x_0\|x^\ast$ as a norming functional of $x_0$. Note that $\|z^\ast\|=\|x_0\|$ and $z^\ast(x_0)=\|x_0\|^2$.
\end{dfn}
If $X$ is a normed space, then by the Hahn-Banach theorem, every $x_0\in X$ has a norming functional.

\begin{dfn}[\textcolor{OliveGreen}{Fr{\'e}chet-G{\^a}teaux derivative}] A map of normed spaces $F:X\ra Y$ is (Fr{\'e}chet-) differentiable at $x\in X$ if there exists a linear map $dF_x:X\ra Y$ and a continuous map $o_x\in C(X,Y)$ such that
\begin{equation*}
\textstyle F(x+h)=F(x)+dF_xh+o_x(h)~~~~\txt{for all}~~h\in X,~~~~\txt{with}~~~~\lim\limits_{\|h\|\ra0}{\|o_x(h)\|\over\|h\|}=0.
\end{equation*}
The map $dF:X\ra L(X,Y)$, $x\mapsto  dF_x$ is called the (Fr{\'e}chet) derivative of $F$, and the linear map $dF_x:X\ra Y$ is called the (Fr{\'e}chet) derivative of $F$ at $x$.

    When the limit is required to exist only ``linearly'' (i.e., in one direction at a time), we get a weaker (G{\^a}teaux) version of the derivative: $F$ is G{\^a}teaux-differentiable at $x\in X$ if there exists a linear map $DF_x:X\ra Y$ and a continuous map $o_x\in C(X,Y)$ such that for every $h\in X$ with $\|h\|=1$,
\begin{equation*}
\textstyle F(x+th)=F(x)+~t~DF_xh+o_x(th),~~~~\txt{for all}~~t\in \Real,~~~~\txt{with}~~\lim\limits_{t\ra0}{\|o_x(th)\|\over|t|}=0.
\end{equation*}
The map $DF:X\ra L(X,Y)$, $x\mapsto DF_x$ is called the G{\^a}teaux derivative of $F$, and the map $D_hF:X\ra Y$, $x\mapsto DF_xh$ is called the directional derivative of $F$ along $h$. Accordingly, the linear map $DF_x:X\ra Y$ is called the G{\^a}teaux derivative of $F$ at $x$, and the vector $DF_xh\in Y$ is called the directional derivative of $F$ at $x$ along $h$.
\end{dfn}

\begin{rmk*}[\textcolor{OliveGreen}{Mean value theorem: Theorem 1.8 on page 13 of \cite{ambro-prodi1993}}] If $F:O\subseteq X\ra Y$ is G{\^a}teaux-differentiable and $O$ is open, then for any $x_1,x_2\in X$ such that ~$[x_1,x_2]:=\big\{\eta(t)=(1-t)x_1+tx_2:t\in[0,1]\big\}\subseteq O$,~
\begin{equation*}
\textstyle\|F(x_1)-F(x_2)\|\leq C(x_1,x_2)\|x_1-x_2\|,~~\txt{where}~~C(x_1,x_2):=\sup\limits_{x\in[x_1,x_2]}\left\|DF_x\right\|.
\end{equation*}
\end{rmk*}

\begin{dfn}[\textcolor{OliveGreen}{Semi-inner products on $X$}]\label{SemiIPs} Let $X$ be a normed space and $x,y\in X$. We define
\begin{equation*}
\textstyle\langle x,y\rangle_-:=\inf\limits_{y^\ast\in\F y}\langle x,y^\ast\rangle,~~~~\langle x,y\rangle_+:=\sup\limits_{y^\ast\in\F y}\langle x,y^\ast\rangle,
\end{equation*}
where $\F:X\ra\P(X^\ast)$ is the duality map of $X$, given by the set of norming functionals
\begin{equation*}
\F x~:=~\left\{x^\ast\in X^\ast:\|x^\ast\|=\|x\|,~x^\ast(x)=\|x\|^2\right\},~~~~\txt{for all}~~x\in X.
\end{equation*}
\end{dfn}

\begin{prp}[\textcolor{OliveGreen}{Derivative characterization of semi-inner products}]\label{SemiInnCh}
Let $X$ be a normed space. The semi-inner products in Definition \ref{SemiIPs} are determined by one-sided derivatives of the norm as follows.
\begin{equation}
\label{SemiInnEq}\textstyle\langle x,y\rangle_-=\|y\|\lim\limits_{t\uparrow 0}{\|y+tx\|-\|y\|\over t},~~~~\langle x,y\rangle_+=\|y\|\lim\limits_{t\downarrow 0}{\|y+tx\|-\|y\|\over t}.
\end{equation}
\end{prp}
\begin{proof}
See \cite{deimling}, Proposition 12.3(d), page 115.
\end{proof}

\begin{lmm}[\textcolor{OliveGreen}{Derivative of the norm along trajectories}]\label{AbsCont}
If $X$ is a normed space and $\gamma:[0,1]\ra X$ is a $C^1$-smooth path, then the following are true.
\bit
\item[(i)] The function $\vphi:[0,1]\ra\Real$, $\vphi(t)=\|\gamma(t)\|$ is absolutely continuous, i.e.,
\begin{equation}
\label{AbsContEq}\textstyle\vphi' ~~~~\txt{exists a.e.},~~~~\vphi'\in L([0,1]), ~~~~\txt{and}~~~~ \vphi(t)=\vphi(0)+\int_0^t\vphi'(s)ds.
\end{equation}
\item[(ii)] With $\gamma'(t)={d\over dt}\gamma(t):=d\gamma_t$, the derivative of $\vphi$ can be expressed in the following form: For a.e. $t\in[0,1]$,
\begin{equation}
\label{LinTrans}\textstyle\vphi'(t)=\lim\limits_{h\downarrow 0}{\|\gamma(t)+h\gamma'(t)\|-\|\gamma(t)\|\over h}~\sr{(\ref{SemiInnEq})}{=}~{1\over\|\gamma(t)\|}\big\langle\gamma'(t),\gamma(t)\big\rangle_+.
\end{equation}
\eit
\end{lmm}
\begin{proof}
(i) Let $C:=\sup_{[0,1]}\|\gamma'\|$, where $\gamma'$ is the derivative of $\gamma$. Then the mean value theorem gives
\begin{equation*}
\big|\vphi(a)-\vphi(b)\big|=\big|\|\gamma(a)\|-\|\gamma(b)\|\big|\leq\big\|\gamma(a)-\gamma(b)\big\|\leq C|a-b|,
\end{equation*}
which shows $\vphi$ is absolutely continuous.

(ii) If $\vphi$ is differentiable at $t\in[0,1]$, then
\begin{equation*}
\textstyle {d\over dt}\|\gamma(t)\|=\lim\limits_{h\ra0}{\|\gamma(t+h)\|-\|\gamma(t)\|\over h}~\sr{(s)}{=}~\lim\limits_{h\ra 0}{\|\gamma(t)+h\gamma'(t)\|-\|\gamma(t)\|\over h},
\end{equation*}
where step (s) holds because $\gamma(t+h)=\gamma(t)+h\gamma'(t)+o_t(h)$, and so
\begin{equation*}
\textstyle \lim\limits_{h\ra 0}\left|{\|\gamma(t+h)\|-\|\gamma(t)+h\gamma'(t)\|\over h}\right|\leq \lim\limits_{h\ra 0}{\|o_t(h)\|\over|h|}=0.\qedhere
\end{equation*}
\end{proof}

\begin{lmm}[\textcolor{OliveGreen}{Semi-monotonicity of the radial projection}]\label{accret-lmm}
If $X$ is a normed space, the map $X\backslash\{0\}\ra X$ given by $x\mapsto \hat{x}:={x\over\|x\|}$ satisfies
\begin{equation}
\label{AccretEq}\langle\hat{x}-\hat{y},x-y\rangle_-\geq0~~~~\txt{for all}~~~~x,y\in X\backslash\{0\}.
\end{equation}
\end{lmm}
\begin{proof}
Fix any two vectors $x,y\in X\backslash\{0\}$. If $\|x\|=\|y\|$, then $\langle\hat{x}-\hat{y},x-y\rangle_-=\|x\|^{-1}\langle x-y,x-y\rangle_-\geq0$. So, assume $\|x\|>\|y\|$. Consider the convex function
\begin{equation*}
\vphi(t)=\|x-y+t(\hat{x}-\hat{y})\|=\big\|(t+\|x\|)\hat{x}-(t+\|y\|)\hat{y}\big\|,~~~~\txt{for}~~t\in\Real.
\end{equation*}
By Proposition \ref{SemiInnCh}, it is enough to show that the left-sided derivative of $\vphi$ is nonnegative at $t=0$.

Observe that $\vphi(-\|x\|)=\vphi(-\|y\|)= \|x\|-\|y\|$. Since $\vphi$ is convex and $\|x\|\neq \|y\|$, it follows that $\vphi$  attains its minimum on $\big[-\|x\|, -\|y\|\big]$. So, $\vphi$ is nondecreasing on $\big[-\|y\|,+\infty\big)$. Hence, both one-sided derivatives of $\vphi$ are nonnegative at $t=0$.
\end{proof}

\begin{thm}[\textcolor{OliveGreen}{Analog of Theorem 1.1 in \cite{kovalev2016}}]\label{MainThm}
If $X$ is a normed space, then for each $n\geq 2$ there exists a retraction ~$r:X(n)\ra X(n-1)$ that is H{\"o}lder-continuous on bounded subsets of $X(n)$.
\end{thm}
\begin{proof}
We will proceed in six steps to construct the retraction and prove its continuity.

{\flushleft \ul{1. {\bf Evolution equation and collision time}}}:~ Equip $X^n$ with the metric {\small $d\big(x,y\big)=\sum_{i=1}^n\|x_i-y_i\|$}, which makes $X^n$ a normed space. Let $D:=\{x\in X^n:x_i=x_j~\txt{for some}~i\neq j\}$, and consider the map
\begin{align*}
J&\textstyle=(J_1,...,J_n):X^n\backslash D\ra X^n,~~~~~~~~J_i(x):=\sum_{j\neq i}{x_i-x_j\over\|x_i-x_j\|},\\
\|J\|&\textstyle=\sum_i\|J_i\|=\sum_i\left\|\sum_{j\neq i}{x_i-x_j\over\|x_i-x_j\|}\right\|\leq n(n-1).
\end{align*}
Note that ~$X(n)\backslash X(n-1)=\big\{x=\{x_1,...,x_n\}:~(x_1,...,x_n)\in X^n\backslash D\big\}$. Consider the system of ordinary differential equations
\begin{align}
\label{KOVeq2}\textstyle{du_i(t)\over dt}&=-J_i\big(u(t)\big),~~~~u_i(0)=x_i,~~~~i=1,...,n,\\
\label{KOVeq3}\textstyle\left\|{du_i(t)\over dt}\right\|&=\|J_i\big(u(t)\big)\|\leq n-1,~~~~i=1,...,n.
\end{align}
Beginning with $u(0)=x\in X^n\backslash D$, by Picard's theorem, the system (\ref{KOVeq2}) continues to have a unique solution
\begin{equation*}
u(t)\in X^n\backslash D,~~~~\txt{with each}~~~~u_i(t)\in\txt{Span}\{x_1,...,x_n\}~~\txt{in}~~X,
\end{equation*}
until we reach the set $D$ (a situation we will refer to as ``\emph{collision}''). Let ~$T(x):=\sup\{t:~t\geq 0,~u(t)\in X^n\backslash D\}$,~ i.e., $[0,T(x))$ is the maximal interval of existence of the solution of (\ref{KOVeq2}).
\begin{rmk}\label{CommColl}
Note that for any $0<\tau<T(x)$, $u(t+\tau)$ is the unique solution of the system
\begin{equation*}
\textstyle {dw_i(t)\over dt}=-J_i\big(w(t)\big),~~~~w_i(0)=u_i\big(\tau\big).
\end{equation*}
This implies the point of collision for $w$ is the same as for $u$, i.e., {\small $u\big(T(x)\big)=w\Big(T\big(u(\tau)\big)\Big)=u\Big(T\big(u(\tau)\big)+\tau\Big)$}. Equivalently, we have
\begin{equation}
\label{TimeTrans}T(x)=T\big(u(\tau)\big)+\tau,~~\txt{or}~~T\big(u(\tau)\big)=T(x)-\tau,~~~~\txt{for all}~~~~0<\tau<T(x).
\end{equation}
\end{rmk}

{\flushleft \ul{2. {\bf Bounds on the collision time $T(x)$}}}:~ With $\delta$ as in Definition \ref{MinSep},
\begin{align*}
\delta(x)&=\min_{i\neq j}\|x_i-x_j\|\leq \min_{i\neq j}\Big(\|x_i-u_i(T(x))\|+\|u_i(T(x))-u_j(T(x))\|+\|u_j(T(x))-x_j\|\Big)\\
  &\sr{(s)}{\leq} (n-1)T(x)+\min_{i\neq j}\|u_i(T(x))-u_j(T(x))\|+(n-1)T(x)\\
  &=(n-1)T(x)+0+(n-1)T(x)=2(n-1)T(x),
\end{align*}
where the mean value theorem is used at step (s). Therefore,
\begin{equation*}
\textstyle T(x)\geq{\delta(x)\over 2(n-1)}.
\end{equation*}

Renumbering the points $x_i$, we may assume $\delta(x)=\|x_1-x_2\|$. Let $\vphi(t):=\|u_1(t)-u_2(t)\|$. Then by Lemma \ref{AbsCont}, $\vphi$ is absolutely continuous and, for all $t$, its derivative satisfies
\begin{align*}
\vphi'(t)~&\textstyle\sr{(\ref{LinTrans})}{=}\left\langle{du_1\over dt}-{du_2\over dt},{u_1-u_2\over\|u_1-u_2\|}\right\rangle_+\sr{(\ref{KOVeq2})}{=}
-\left\langle J_1(u)-J_2(u),{u_1-u_2\over \|u_1-u_2\|}\right\rangle_-\nn\\
   &\textstyle=-\left\langle\sum_{j\neq 1}{u_1-u_j\over\|u_1-u_j\|}-\sum_{j\neq 2}{u_2-u_j\over\|u_2-u_j\|},{u_1-u_2\over\|u_1-u_2\|}\right\rangle_-\nn\\
   &\textstyle~=-2-\sum_{j= 3}^n{\left\langle{u_1-u_j\over\|u_1-u_j\|}-{u_2-u_j\over\|u_2-u_j\|},(u_1-u_j)-(u_2-u_j)\right\rangle_-\over\|u_1-u_2\|}\nn\\
   &\textstyle\sr{(\ref{AccretEq})}{\leq} -2,~~~~\txt{for almost all}~~~~0<t<T(x),
\end{align*}
where step (\ref{AccretEq}) refers to the property of the radial projection proved in Lemma \ref{accret-lmm}. Upon integration of the above inequality, we get $\vphi\big(T(x)\big)-\vphi(0)\leq -2T(x)$, which implies
\begin{equation}
\label{KOVeq4}\textstyle T(x)~\leq~{\delta(x)\over 2}.
\end{equation}

{\flushleft \ul{3. {\bf Definition of the retraction}}}:~ Define $r:X(n)\ra X(n-1)$ as follows. If $x\in X(n)\backslash X(n-1)$, let $r(x)=r\big(\{x_i\}\big):=\left\{u_i\big(T(x)\big)\right\}=u\big(T(x)\big)$, which is a well defined map since the order of enumeration is unimportant. If $x\in X(n-1)$, let $r(x):=x$. Then $r|_{X(n-1)}=id_{X(n-1)}$. It remains to show that $r$ is continuous. Specifically, we will show that for all $x,y\in X(n)$,
\begin{equation}
\label{KOVeq5}d_H\big(r(x),r(y)\big)~\leq~n(2n-1)\diam\left(x\cup y\right)^{1-{1\over 2n-1}}d_H\left(x,y\right)^{1\over 2n-1}.
\end{equation}
For $x,y\in X^n$, let $u,v$ be the solutions of (\ref{KOVeq2}) with initial data $u(0)=x$, $v(0)=y$. Recall that for all $x,y\in X(n)$,
\begin{equation*}
d_H\big(x,y\big)=\max\left\{\sup_{1\leq i\leq n}\inf_{1\leq j\leq n}\|x_i-y_j\|,\sup_{1\leq i\leq n}\inf_{1\leq j\leq n}\|x_j-y_i\|\right\}.
\end{equation*}
{\flushleft \ul{4. {\bf Estimate of Hausdorff distance to the collision point}}}:~ Using (\ref{KOVeq3}) and (\ref{KOVeq4}),
\begin{align*}
d_H&\big(r(x),x\big)=\max\left\{\sup_i\inf_j\|u_i(T(x))-x_j\|,\sup_i\inf_j\|u_j(T(x))-x_i\|\right\}\nn\\
   \leq&\max\left\{\sup_i\inf_j\Big((n-1)T(x)+\|x_i-x_j\|\Big),\sup_i\inf_j\Big((n-1)T(x)+\|x_j-x_i\|\Big)\right\}\nn\\
   =&\max\left\{\Big((n-1)T(x)+\sup_i\inf_j\|x_i-x_j\|\Big),\Big((n-1)T(x)+\sup_i\inf_j\|x_j-x_i\|\Big)\right\}\nn\\
   =&\textstyle\max\left\{\Big((n-1)T(x)+0\Big),\Big((n-1)T(x)+0\Big)\right\}=(n-1)T(x)\leq(n-1){\delta(x)\over 2},
\end{align*}
and a similar bound holds for $d(r(y),y)$. Therefore,

\begin{equation}
\label{KOVeq6}\textstyle d_H\big(r(x),x\big)\leq{n-1\over 2}\delta(x),~~~~d_H\big(r(y),y\big)\leq{n-1\over 2}\delta(y).
\end{equation}

{\flushleft \ul{5. {\bf Estimate of Hausdorff distance at first collision}}}: By Lemma \ref{AbsCont}, the function $g(t):=\sum_{i=1}^n\|u_i(t)-v_i(t)\|$ is absolutely continuous and its derivative satisfies the following: For a.e. $0<t<T:=\min\Big(T(x),T(y)\Big)$,
{\small\begin{align*}
g'(t)\sr{(\ref{LinTrans})}{=}\textstyle &\textstyle\sum\limits_i\left\langle{du_i\over dt}-{dv_i\over dt},{u_i-v_i\over\|u_i-v_i\|}\right\rangle_+\sr{(\ref{KOVeq2})}{=}-\sum\limits_i\left\langle J_i(u)-J_i(v),{u_i-v_i\over\|u_i-v_i\|}\right\rangle_-,\\
|g'(t)|\leq&\textstyle \sum\limits_i\|J_i(u)-J_i(v)\|\leq \sum\limits_i\sum\limits_{j\neq i}\left\|{u_i-u_j\over\|u_i-u_j\|}-{v_i-v_j\over\|v_i-v_j\|}\right\|= 2\sum\limits_{i<j}\left\|{u_i-u_j\over\|u_i-u_j\|}-{v_i-v_j\over\|v_i-v_j\|}\right\|\\
\sr{(s)}{\leq}&\textstyle 4{\sum_{i<j}\left\|(u_i-u_j)-(v_i-v_j)\right\|\over\max(\|u_i-u_j\|,\|v_i-v_j\|)}
   \leq {4(n-1)\sum_i\|u_i-v_i\|\over\min\limits_{i<j}\Big(\max(\|u_i-u_j\|,\|v_i-v_j\|)\Big)}
   \leq {4(n-1)\over\max\Big(\delta\big(u(t)\big),\delta\big(v(t)\big)\Big)}g(t)\\
   \sr{(\ref{KOVeq4})}{\leq}&\textstyle {2(n-1)\over\max\Big(T\big(u(t)\big),T\big(v(t)\big)\Big)}g(t)
   \sr{(\ref{TimeTrans})}{=}{2(n-1)\over\max\Big(T(x),T(y)\Big)-t}g(t),\nn
\end{align*}}
where step (s) is due to the inequality $\left|{x\over\|x\|}-{y\over\|y\|}\right|\leq{2\|x-y\|\over\max(\|x\|,\|y\|)}$ from \cite{dunkl1964,thele1974}. It follows that
\begin{equation*}
\textstyle |g'(t)|\leq{2(n-1)\over T-t}g(t),~~~~\txt{for a.e.}~~0<t<T.
\end{equation*}
Since we also have ~$|g'(t)|\leq 2(n-1)$~ for a.e. $0<t<T$, it follows that ~$|g'(t)|\leq 2(n-1)\min\left({1\over T-t}g(t),1\right)$~ for a.e. $0<t<T$. We consider various cases as follows.
\bit
\item[(i)]\ul{$g(\tau)=T-\tau$ for some $\tau\in(0,T)$}:~ The bound $g'(t)\leq{2(n-1)\over T-t}g(t)$ for a.e. implies
\begin{equation*}
\textstyle g(\tau)\leq \exp\left(\int_0^\tau{2(n-1)\over T-s}ds\right)g(0)=
\textstyle \left(\frac{T}{T-\tau}\right)^{2(n-1)} g(0)=\left(\frac{T}{g(\tau)}\right)^{2(n-1)} g(0),
\end{equation*}
which in turn implies
\begin{equation}
\label{gBound1} g(\tau)\leq T^{1-{1\over 2n-1}}g(0)^{1\over 2n-1}.
\end{equation}
Also, the bound $g'(t)\leq 2(n-1)$ a.e. implies ~$g(T)-g(\tau)\leq 2(n-1)(T-\tau) = 2(n-1)g(\tau)$,~ i.e.,
\begin{equation}
\label{gBound2}g(T)\leq (2n-1)g(\tau)\sr{(\ref{gBound1})}{\leq}(2n-1)T^{1-{1\over 2n-1}}g(0)^{1\over 2n-1}.
\end{equation}
\item[(ii)]\ul{$g(t)<T-t$ for all $t\in (0,T)$}:~ The bound $g'(t)\leq{2(n-1)\over T-t}g(t)$ a.e. shows (\ref{gBound1}) holds for all $t$, that is, for all $0<t<T$,
\begin{equation}
\label{gBound3}g(t)\leq T^{1-{1\over 2n-1}}g(0)^{1\over 2n-1}\leq (2n-1)T^{1-{1\over 2n-1}}g(0)^{1\over 2n-1}.
\end{equation}
\item[(iii)]\ul{$g(t)> T-t$ for all $t\in(0,T)$}:~  The bound $g'(t)\leq 2(n-1)$ a.e. implies ~$g(T)-g(t)\leq 2(n-1)(T-t)< 2(n-1)g(t)$,~ i.e., $g(T)< (2n-1)g(t)$, for all $0<t<T$. Therefore, in the limit $t\ra0$, we get
\begin{equation}
\label{gBound4}g(T)\leq (2n-1)g(0).
\end{equation}
\eit
Since $T\leq{\min\big(\delta(x),\delta(y)\big)\over 2}<\diam(x\cup y)$ and $\diam(x\cup y)\geq g(0)$, it follows that all three cases above imply
\begin{equation}
\label{gBound5}g(T)\leq(2n-1)~\diam(x\cup y)^{1-{1\over 2n-1}}~g(0)^{1\over 2n-1}~=~C_n(x,y)~g(0)^{1\over 2n-1},
\end{equation}
where $C_n(x,y):=(2n-1)~\diam(x\cup y)^{1-{1\over 2n-1}}$. Therefore,
\begin{align*}
d_H\big(u(T)&,v(T)\big)\leq d\big(u(T),v(T)\big)=\max_i\|u_i(T)-v_i(T)\|\sr{(\ref{gBound5})}{\leq}~C_n(x,y)\max_i\|u_i(0)-v_i(0)\|^{1\over 2n-1},
\end{align*}
which implies
\begin{equation}
\label{KOVeq7} d_H\big(u(T),v(T)\big)\leq~C_n(x,y)\max_i\|x_i-y_i\|^{1\over 2n-1}.
\end{equation}

{\flushleft \ul{6. {\bf Estimate of Hausdorff distance between collision points}}}: We consider two cases as follows, where $\rho:=d_H(x,y)$.
\bit
\item \ul{Case 1}: $\delta(x)+\delta(y)\leq 4\rho$. In this case, we obtain (\ref{KOVeq5}) as follows.
\begin{align*}
d_H\big(&r(x),r(y)\big)\textstyle\leq d_H\big(r(x),x\big)+d_H(x,y)+d_H(y,r(y))\sr{(\ref{KOVeq6})}{\leq}~ {n-1\over 2}\delta(x)+d_H(x,y)+{n-1\over 2}\delta(y)\\
  &\leq \rho+2(n-1)\rho=(2n-1)\rho\leq n(2n-1)\diam\left(x\cup y\right)^{1-{1\over 2n-1}}d_H\left(x,y\right)^{1\over 2n-1}.
\end{align*}
\item \ul{Case 2}: $\delta(x)+\delta(y)>4\rho$. In this case, $\delta(x)>2\rho$ or $\delta(y)>2\rho$, and so (\ref{KOVeq8}) applies. By definition, $r(x)=u(T(x))$. Let $z:=v(T(x))$. Then $d_H\big(r(x),z\big)=d_H\Big(u\big(T(x)\big),v\big(T(x)\big)\Big)\sr{(\ref{KOVeq7})}{\leq}C_n(x,y)\max_i\|x_i-y_i\|^{1\over 2n-1}\sr{(\ref{KOVeq8})}{\leq} C_n(x,y)\rho^{1\over 2n-1}$. That is,
\begin{align}
\label{KOVeq9}d_H\big(r(x),z\big)\leq C_n(x,y)\rho^{1\over 2n-1}.
\end{align}
Since $\delta$ is 2-Lipschitz (Lemma \ref{DeltaCont}) and $\delta\big(r(x)\big)=0$, with $z_i:=v_i\big(T(x)\big)$,
\begin{equation*}
\delta(z)=|\delta\big(r(x)\big)-\delta(z)|\leq 2d_H\big(r(x),z\big)\sr{(\ref{KOVeq9})}{\leq}2C_n(x,y)\rho^{1\over 2n-1},
\end{equation*}
which together with (\ref{KOVeq6}) implies
\begin{equation}
\label{KOVeq10}\textstyle d_H\big(r(z),z\big)\sr{(\ref{KOVeq6})}{\leq}{n-1\over 2}\delta(z)\leq (n-1)C_n(x,y)\rho^{1\over 2n-1}.
\end{equation}
Using Remark \ref{CommColl} (which says $r(z)=r(y)$) at step (R\ref{CommColl}) below,
\begin{align*}
d_H\big(r(x),r(y)\big)&\sr{\txt{(R\ref{CommColl})}}{=}d_H\big(r(x),r(z)\big)\leq d_H\big(r(x),z\big)+d_H\big(r(z),z\big)\sr{(\ref{KOVeq9}),(\ref{KOVeq10})}{\leq} n C_n(x,y)\rho^{1\over 2n-1},
\end{align*}
which in turn implies (\ref{KOVeq5}), i.e.,
\begin{equation*}
d_H\big(r(x),r(y)\big)\leq n(2n-1)\diam\left(x\cup y\right)^{1-{1\over 2n-1}}d_H\left(x,y\right)^{1\over 2n-1}.
\end{equation*}
\eit
This completes the proof.
\end{proof}

\begin{rmk}[\textcolor{OliveGreen}{Connection with \cite{kovalev2016}}]
When $X$ is a Hilbert space, the norm is differentiable. Thus, with the metric $d(x,y)=\left(\sum_i\|x_i-y_i\|^2\right)^{1/2}$ on $X^n$, the function $g(t)={1\over 2}\sum_i\|u_i(t)-v_i(t)\|^2$ satisfies the following (for all $0<t<T$):
{\small\begin{align*}
g'(t)&\textstyle =\sum\limits_i\left\langle{du_i\over dt}-{dv_i\over dt},u_i-v_i\right\rangle_+=-\sum_i\Big\langle J_i(u)-J_i(v),u_i-v_i\Big\rangle_-\\
  &\textstyle =-\sum\limits_{1\leq i<j\leq n}\left\langle{u_i-u_j\over\|u_i-u_j\|}-{v_i-v_j\over\|v_i-v_j\|},(u_i-v_i)^\ast-(u_j-v_j)^\ast\right\rangle\\
  &\textstyle \sr{(a)}{=}-\sum\limits_{1\leq i<j\leq n}\left\langle{u_i-u_j\over\|u_i-u_j\|}-{v_i-v_j\over\|v_i-v_j\|},(u_i-u_j)^\ast-(v_i-v_j)^\ast\right\rangle\sr{(b)}{\leq}0,
\end{align*}}where step (a) is due to linearity of the duality map (by Riesz representation theorem), and step (b) is due to monotonicity of the radial projection $F(x)=x/\|x\|$, $x\neq0$. This leads to $g(t)\leq g(0)$ for all $0<t<T$. Hence, as in \cite{kovalev2016}, the map $r:X(n)\ra X(n-1)$ is a Lipschitz retraction.
\end{rmk}

\section{Quasiconvexity of Finite Subset Spaces of a Geodesic Metric Space}\label{QuasiConv}

In this section, unless it is stated otherwise, $X$ is any metric space.


\begin{dfn}[\textcolor{OliveGreen}{Spaced Pairs of Points in a Metric Space}]
Let $(X,d)$ be a metric space. Two points $x,y\in X$ are spaced (or form a spaced pair) if
\begin{equation*}
\ol{N}_r(x)\cap\ol{N}_r(y)=\emptyset~~~~\txt{for all}~~~~0<r<d(x,y).
\end{equation*}
\end{dfn}
Equivalently, $x,y\in X$ are spaced $\iff$ $d(x,y)\leq \max\big\{d(x,z),d(z,y)\big\}$ for all $z\in X$.

\begin{lmm}\label{SpacedPairLmm}
If $X$ is a normed space, then for $n\geq 3$ the metric space $X(n)$ contains spaced pairs.
\end{lmm}
\begin{proof}
Since every normed space contains a copy of $\Real$, it is enough to show that $\Real(n)$ contains spaced pairs.

Fix a number $m>3$. Let $x=\{x_1,...,x_n\}\in \Real(n)\backslash \Real(n-1)$ be given by $x_1=0$, $x_2=m-1$, $x_3=m+1$, and $x_i=(i-2)m+1$ for $i\geq 4$. Similarly, let $y=\{y_1,...,y_n\}\in \Real(n)\backslash \Real(n-1)$ be given by $y_1=-1$, $y_2=1$, $y_3=m$, and $y_i=(i-2)m+2$ for $i\geq 4$. That is,
\begin{align*}
x&=\{x_1,...,x_n\}=\{0,m-1,m+1,2m+1,3m+1,\cdots, (n-2)m+1\},\\
y&=\{y_1,...,y_n\}=\{-1,1,m,2m+2,3m+2,\cdots, (n-2)m+2\}.
\end{align*}
The points $x,y$ form a spaced pair because $\ol{N}_r(x)\cap \ol{N}_r(y)=\emptyset$ for all $0<r<1=d_H(x,y)$, where the proof is as follows.

Suppose on the contrary that $z=\{z_1,...,z_n\}\in\ol{N}_r(x)\cap\ol{N}_r(y)$ for some $0<r<1$. Consider the sets $A_1=\{y_1,x_1,y_2\}=\{-1,0,1\}$, $A_2=\{x_2,y_3,x_3\}=\{m-1,m,m+1\}$, and $A_k=\{x_k,y_k\}=\{(k-1)m+1,(k-1)m+2\}$ for $3\leq k\leq n-1$. Then $A_1,A_2$ each contain at least two elements of $z$, while $A_3,...,A_{n-1}$ each contain at least one element of $z$, i.e., at least $2+2+(n-3)=n+1$ elements of $z$ are required (a contradiction).
\end{proof}

\begin{rmk}\label{SpacedPairRmk}
Note that if $X$ is any geodesic space, then Lemma \ref{SpacedPairLmm} remains true. To prove this, let some geodesic segment play the role of the real line in the proof of Lemma \ref{SpacedPairLmm}.
\end{rmk}

\begin{dfn}[\textcolor{OliveGreen}{Geodesic}]
Let $X$ be a metric space. A path $\gamma:[0,1]\ra X$ is a geodesic if $d(\gamma(t),\gamma(t'))=d(\gamma(0),\gamma(1))|t-t'|$ for all $t,t'\in[0,1]$.
\end{dfn}
An equivalent definition of a geodesic in terms of paths that are parameterized by arc length can be found in \cite{papado2014}, Definition 2.2.1 (page 56).
\begin{dfn}[\textcolor{OliveGreen}{Length of a path}]
Let $X$ be a metric space. Given a path $\gamma:[0,1]\ra X$, the length of $\gamma$ is
\begin{equation*}
l(\gamma)~:=~\sup\big\{l_P(\gamma):P\subseteq [0,1]~~\txt{a finite partition}\big\},
\end{equation*}
where $l_P(\gamma):=\sum_{i=1}^kd(\gamma(t_{i-1}),\gamma(t_i))$ is the length of $\gamma$ over $P=\{0=t_0<t_1<\cdots<t_k=1\}$.
\end{dfn}

\begin{lmm}
Let $X$ be a metric space and $\gamma:[0,1]\ra X$ a path. Then (i) $\gamma$ is a geodesic $\iff$ (ii) $\gamma$ is injective and $l(\gamma)=d(\gamma(0),\gamma(1))$ $\iff$ (iii) $d(\gamma(t),\gamma(t'))\leq d(\gamma(0),\gamma(1))|t-t'|$ for all $t,t'\in[0,1]$.
\end{lmm}
\begin{proof}
See \cite{papado2014}, Section 2.2, pages 56-60.
\end{proof}

\begin{dfn}[\textcolor{OliveGreen}{Quasigeodesic, Quasiconvex space, Geodesic space}]
Let $X$ be a metric space and $\ld\geq1$. A path $\gamma:[0,1]\ra X$ is a $\ld$-quasigeodesic if $d(\gamma(t),\gamma(t'))\leq\ld d(\gamma(0),\gamma(1))|t-t'|$ for all $t,t'\in[0,1]$. In particular, $\gamma$ is a geodesic $\iff$ $\gamma$ is a $1$-quasigeodesic.

We say $X$ is a $\ld$-quasiconvex space if for every $x,y\in X$, there exists a $\ld$-quasigeodesic $\gamma:[0,1]\ra X$ from $x$ to $y$, i.e., such that $\gamma(0)=x$, $\gamma(1)=y$. A $1$-quasiconvex space is called a geodesic space.
\end{dfn}
Note that a $\ld$-quasigeodesic is also called a $\ld$-quasiconvex path, \cite{hakobyan-herron2008}, page 205. In \cite{tyson-wu2005}, page 317, a quasigeodesic is differently defined to be a path that is a bi-Lipschitz embedding. Injectivity of the path is not required in our definition.

\begin{dfn}[\textcolor{OliveGreen}{Complete relation, Incomplete relation, Proximal relation}]
Let $A,B$ be sets. A relation on $A$ and $B$ is a set $R\subseteq A\times B$. Given a relation $R\subseteq A\times B$, let
\begin{align*}
A_b(R)&=\{a\in A:(a,b)\in R\},~~~~\txt{($R$-correspondents of $b$ in $A$)}\\
B_a(R)&=\{b\in B:(a,b)\in R\},~~~~\txt{($R$-correspondents of $a$ in $B$)}\\
A(R)&\textstyle =\bigcup_{b\in B}A_b(R)=\{a\in A:(a,b)\in R~\txt{for some}~b\in B\},~~~~\txt{(Left projection of $R$)}\\
B(R)&\textstyle =\bigcup_{a\in A}B_a(R)=\{b\in B:(a,b)\in R~\txt{for some}~a\in A\},~~~~\txt{(Right projection of $R$)}.
\end{align*}
Then, we say $R$ is complete if $A(R)=A$ and $B(R)=B$. Otherwise, $R$ is incomplete.

If $X$ is a metric space and $A,B\subseteq X$, a relation $R\subseteq A\times B$ is proximal if $d(a,b)\leq d_H(A,B)$ for all $(a,b)\in R$.
\end{dfn}

Note that by the definition of Hausdorff distance, for any $x,y\in X(n)$ there exists a proximal complete relation $R\subseteq x\times y$. This knowledge will be used in the proof of Corollary \ref{ProxSplitCrl}.

\begin{dfn}[\textcolor{OliveGreen}{Orders of an element of a relation}]
Let $R\subseteq A\times B$ be relation and $(a,b)\in R$. The left and right orders of $(a,b)$ in $R$ are
\begin{equation*}
O_l(a,b):=|A_b(R)|=\big|(A\times\{b\})\cap R\big|,~~~~O_r(a,b):=|B_a(R)|=\big|(\{a\}\times B)\cap R\big|.
\end{equation*}
\end{dfn}

\begin{dfn}[\textcolor{OliveGreen}{Essential and Inessential elements of a Complete Relation}]
Let $R\subseteq A\times B$ be a complete relation. We say $(a,b)\in R$ is essential if $O_l(a,b)=1$ or $O_r(a,b)=1$. Otherwise, we say $(a,b)$ is inessential.
\end{dfn}
Note that if $R\subseteq A\times B$ is complete, then an element $(a,b)\in R$ is essential (resp. inessential) $\iff$ the relation $R\backslash\{(a,b)\}\subseteq A\times B$ is incomplete (resp. complete).

\begin{dfn}[\textcolor{OliveGreen}{Reduced Complete Relation}]
We say a complete relation $R\subseteq A\times B$ is reduced (or in reduced form) if every element of $R$ is essential.
\end{dfn}

\begin{prp}[\textcolor{OliveGreen}{Characterization of Reduced Complete Relations}]\label{RedRelDecomp}
Let $R\subseteq A\times B$ be a complete relation. Then $R$ is reduced $\iff$ there exist disjoint union decompositions $A=A'\sqcup A_0\sqcup A''$, $B=B'\sqcup B_0\sqcup B''$ and maps $f:A'\ra B'$, $g:B''\ra A''$, $h:A_0\ra B_0$ such that $f,g$ are surjective, $h$ is bijective, and
\begin{equation*}
R=\Big\{\big(a,f(a)\big):a\in A'\Big\}\sqcup \Big\{\big(a,h(a)\big):a\in A_0\Big\}\sqcup\Big\{\big(g(b),b\big):b\in B''\Big\}.\nn
\end{equation*}
\end{prp}
\begin{proof}
$\Ra$: Assume that $R$ is reduced. Define sets $A_1',B_1''$ and maps $f_1:A_1'\ra B$, $g_1:B_1''\ra A$ as follows.
\begin{align*}
A_1'&:=\left\{a\in A:\big|B_a(R)\big|=1\right\},~~~~f_1(a):=\txt{the unique element $b\in B$ such that}~(a,b)\in R.\\
B_1''&:=\left\{b\in B:\big|A_b(R)\big|=1\right\},~~~~g_1(b):=\txt{the unique element $a\in A$ such that}~(a,b)\in R.
\end{align*}
Then, with $R_1:=\Big\{\big(a,f_1(a)\big):a\in A_1'\Big\}$, $R_2:=\Big\{\big(g_1(b),b\big):b\in B_1''\Big\}$, we have
\begin{equation*}
R=R_1\cup R_2=(R_1\backslash R_0)\sqcup R_0\sqcup(R_2\backslash R_0),
\end{equation*}
where $R_0:=R_1\cap R_2$ consists of elements $(u,v)\in R$ such that
\begin{align*}
(u,v)&=\big(a,f_1(a)\big)=\big(g_1(b),b\big)~~\Big[~\iff~~u=a=g_1(b),~~v=b=f_1(a)~\Big]~~\txt{for some}~~a\in A_1',~~b\in B_1''.
\end{align*}
Let $A_0:=A_1'\cap g_1(B_1'')$, $B_0:=f_1(A_1')\cap B_1''$. Then $f_1(A_0)=B_0$, $g_1(B_0)=A_0$, since
\begin{align*}
A_0=&\Big\{a\in A:~a=g_1(b)~\txt{for some}~b\in B,~\txt{with unique}~(a,f_1(a))\in R,~(g_1(b),b)\in R\Big\}\\
   =&\Big\{a\in A:~\txt{there exists}~b\in B~\txt{such that}~a=g_1(b),~b=f_1(a)\Big\}\\
   =&\Big\{a\in A:~\txt{there exists}~b\in B~\txt{such that}~(a,f_1(a))=(g_1(b),b)\in R\Big\},
\end{align*}
and similarly,
\begin{align*}
B_0=\Big\{b\in B:~\txt{there exists}~a\in A~\txt{such that}~(a,f_1(a))=(g_1(b),b)\in R\Big\}.
\end{align*}
It follows that $f_1|_{A_0}:A_0\ra B_0$, $g_1|_{B_0}:B_0\ra A_0$ are mutually inverse bijections, and
\begin{equation*}
R_0=\big\{(a,f_1(a)):a\in A_0\big\}=\big\{(g_1(b),b):b\in B_0\big\}.
\end{equation*}
Hence, we can set $h=f_1|_{A_0}$, $A'=A_1'\backslash A_0$, $f=f_1|_{A'}$, $B''=B_1''\backslash B_0$, $g=g_1|_{B''}$, $B'=f(A')$, $A''=g(B'')$.

$\La$: The converse is clear by the properties of the maps $f,g,h$.
\end{proof}

\begin{crl}[\textcolor{OliveGreen}{Characterization of Reduced Complete Relations}]\label{RedRelDecomp2}
Let $R\subseteq A\times B$ be a complete relation. Then $R$ is reduced $\iff$ there exist disjoint unions $A=A'\sqcup A''$, $B=B'\sqcup B''$ and surjective maps $f:A'\ra B'$, $g:B''\ra A''$ such that $R=\Big\{\big(a,f(a)\big):a\in A'\Big\}\sqcup\Big\{\big(g(b),b\big):b\in B''\Big\}$.
\end{crl}

\begin{lmm}[\textcolor{OliveGreen}{Reduction of a Finite Complete Relation}]\label{RelReductLmm}
Let $X$ be a metric space and $x,y\in X(n)$. Any (proximal) complete relation $R\subseteq x\times y$ contains a (proximal) reduced complete relation $R_{rc}\subseteq x\times y$, which means $R_{rc}\subseteq R$.
\end{lmm}
\begin{proof}
Since $R\subseteq x\times y$ is finite, we can obtain a reduced complete relation $R_{rc}\subseteq R\subseteq x\times y$ by repeatedly excluding inessential elements of $R$.
\end{proof}

\begin{crl}\label{ProxSplitCrl}
Let $X$ be a metric space and $x,y\in X(n)$. There exist proximal maps $f:x'\subseteq x\ra y$ and $g:y''\subseteq y\ra x$ such that $x=x'\sqcup g(y'')$ and $y=f(x')\sqcup y''$.
\end{crl}
\begin{proof}
By Lemmas \ref{RelReductLmm} and the definition of Hausdorff distance, a proximal reduced  complete relation $R\subseteq x\times y$ exists. Hence, by Corollary \ref{RedRelDecomp2}, the desired proximal maps exist.
\end{proof}

\begin{prp}[\textcolor{OliveGreen}{Sufficient condition for quasigeodesics in $X(n)$}]\label{QGeodExistSuff}
Let $X$ be a geodesic space and $x,y\in X(n)$. If there exists a complete relation $R\subseteq x\times y$  satisfying
\begin{equation}
\label{QGeodCond}|R|\leq n,~~~~d(a,b)\leq \ld d_H(x,y)~~\txt{for all}~~(a,b)\in R,
\end{equation}
then $x,y$ are connected by a $\ld$-quasigeodesic in $X(n)$.

\end{prp}
\begin{proof}
Assume some $R\subseteq x\times y$ satisfies (\ref{QGeodCond}). Then the map $\gamma:[0,1]\ra X(n)$ given by
\begin{equation*}
\gamma(t):=\left\{\gamma_{(a,b)}(t):(a,b)\in R,~\gamma_{(a,b)}~\txt{a geodesic from $a$ to $b$}\right\}\nn
\end{equation*}
is a $\ld$-quasigeodesic from $x$ to $y$, since $\gamma(0)=x$, $\gamma(1)=y$, and for all $t,t'\in[0,1]$ we have
\begin{align*}
d_H(\gamma(t)&,\gamma(t'))=\max\left\{\!\max_{(a,b)\in R}\min_{(c,d)\in R}d\left(\gamma_{(a,b)}(t),\gamma_{(c,d)}(t')\right),\max_{(c,d)\in R}\min_{(a,b)\in R}d\left(\gamma_{(a,b)}(t),\gamma_{(c,d)}(t')\right)\!\right\}\\
   \leq&\max_{(a,b)\in R}d\left(\gamma_{(a,b)}(t),\gamma_{(a,b)}(t')\right)=|t-t'|\max_{(a,b)\in R}d(a,b)\leq\ld|t-t'|d_H(x,y).\qedhere
\end{align*}
\end{proof}

\begin{lmm}[\textcolor{OliveGreen}{Geodesics via proximal reduced complete relations}]\label{GeodPRCRLmm}
Let $X$ be a geodesic space. Then any two finite sets $x,y\subseteq X$ are connected by a geodesic  $\gamma:[0,1]\ra X(N)$, where $N:=\max(|x|,|y|,|x|+|y|-2)$. In particular, any two points $x,y\in X(n)$ are connected by a geodesic in $X\big(\max(n,2n-2)\big)$.
\end{lmm}
\begin{proof}
By Corollaries \ref{RedRelDecomp2} and \ref{ProxSplitCrl}, there exist proximal maps $f:x'\subseteq x\ra y$, $g:y''\subseteq y\ra x$ such that $x=x'\sqcup g(y'')$, $y=f(x')\sqcup y''$, and a proximal reduced complete relation $R\subseteq x\times y$ such that $R=\{(a,f(a)):a\in x'\}\sqcup\{(g(b),b):b\in y''\}$. Thus, we have the following 3 cases. (i) If $x'=\emptyset$, then $|R|=|y|$. (ii) If $y''=\emptyset$,  then $|R|=|x|$. (iii) If $x'\neq\emptyset$, $y''\neq\emptyset$, then
\begin{equation*}
|R|=|x'|+|y''|=|x|+|y|-\big(|f(x')|+|g(y'')|\big)\leq |x|+|y|-2.\nn
\end{equation*}
The conclusion now follows from Proposition \ref{QGeodExistSuff}.
\end{proof}

\begin{thm}[\textcolor{OliveGreen}{Quasiconvexity of $X(n)$}]\label{QConvThm}
If $X$ is a geodesic space, then $X(n)$ is 2-quasiconvex. Moreover, $X(2)$ is a geodesic space, and for $n\geq 3$, $\ld=2$ is the smallest quasiconvexity constant for $X(n)$.
\end{thm}
\begin{proof}
Let $x,y\in X(n)$. If $x=y$, then $x,y$ are connected by the constant path. So, assume $x\neq y$. By Corollary \ref{ProxSplitCrl}, we have proximal maps $f:x'\subseteq x\ra y$, $g:y''\subseteq y\ra x$ such that, with $x'':=g(y'')$ and $y':=f(x')$, we have
\begin{equation*}
x=x'\sqcup g(y'')=x'\sqcup x'',~~~~y=f(x')\sqcup y''=y'\sqcup y''.
\end{equation*}
Let $z:=x''\cup y'$, for which we can verify that $d_H(x,z)\leq d_H(x',y')\leq d_H(x,y)$ and $d_H(z,y)\leq d_H(x'',y'')\leq d_H(x,y)$. Then $R_1=\{(a,f(a)):a\in x'\}\cup\{(c,c):c\in x''\}\subseteq x\times z$ and $R_2=\{(c,c):c\in y'\}\cup\{(g(b),b):b\in y''\}\subseteq z\times y$ are complete relations with the following properties.
\bit
\item $|R_1|\leq n$ and $d(u,v)\leq d_H(x,y)$ for each $(u,v)\in R_1$.
\item $|R_2|\leq n$ and $d(u,v)\leq d_H(x,y)$ for each $(u,v)\in R_2$.
\eit
If $z=x$ or $z=y$, then by Proposition \ref{QGeodExistSuff}, $x,y$ are connected by a quasigeodesic. So, assume $z\neq x$ and $z\neq y$. Then it follows again by Proposition \ref{QGeodExistSuff} that there exists a ${d_H(x,y)\over d_H(x,z)}$-quasigeodesic $\gamma_1:[0,1]\ra X(n)$ from $x$ to $z$, and there exists a ${d_H(x,y)\over d_H(z,y)}$-quasigeodesic $\gamma_2:[0,1]\ra X(n)$ from $z$ to $y$. Let $\gamma=\gamma_1\cdot\gamma_2:[0,1]\ra X(n)$ be the path from $x$ to $y$ given by
\begin{equation*}
\gamma(t):=\left\{
             \begin{array}{ll}
               \gamma_1(2t), & \txt{if}~~t\in[0,1/2] \\
               \gamma_2(2t-1), & \txt{if}~~t\in[1/2,1]
             \end{array}
           \right\}.
\end{equation*}
Then
\begin{equation*}
\textstyle l(\gamma)=l(\gamma_1)+l(\gamma_2)\leq {d_H(x,y)\over d_H(x,z)}d_H(x,z)+{d_H(x,y)\over d_H(z,y)}d_H(z,y)=2d_H(x,y).
\end{equation*}
This shows that $X(n)$ is $2$-quasiconvex.

It follows from Lemma \ref{GeodPRCRLmm} that $X(2)$ is a geodesic metric space. Now let $n\geq 3$. To show $\ld=2$ is the smallest quasiconvexity constant for $X(n)$, let $x,y\in X(n)$ be a spaced pair (which exists by Remark \ref{SpacedPairRmk}). Then $d_H(x,y)\leq\max\left\{d_H(x,z),d_H(z,y)\right\}$ for all $z\in X(n)$. If $\gamma:[0,1]\ra X(n)$ is a $\ld$-quasigeodesic from $x$ to $y$, i.e., $d_H(\gamma(t),\gamma(t'))\leq\ld|t-t'|d_H(x,y)$, then
\begin{align*}
d_H(x,&\gamma(1/2))\textstyle\leq{\ld\over 2}d_H(x,y),~~~~d_H(\gamma(1/2),y)\leq{\ld\over 2}d_H(x,y)\\
~~\Ra&\textstyle~~d_H(x,y)\leq\max\left\{d_H(x,\gamma(1/2)),d_H(\gamma(1/2),y)\right\}\leq{\ld\over 2}d_H(x,y),\\
~~\Ra&~~\ld\geq 2.\qedhere
\end{align*}
\end{proof}

\section*{Acknowledgements}
I am indebted to Leonid Kovalev for important discussions and critical suggestions throughout the preparation of this manuscript. I am also grateful to the anonymous referee whose careful review, with several detailed comments and suggestions, greatly improved readability of the paper.




\end{document}